\documentclass[reqno]{amsart}
\pdfoutput=1
\usepackage{amssymb}
\usepackage{graphicx}
\usepackage[all]{xy}

\begin{document}
\newcommand*{\threesim}{\mathrel{\vcenter{\offinterlineskip\hbox{$\sim$}\vskip-.35ex\hbox{$\sim$}\vskip-.35ex\hbox{$\sim$}}}}

\newtheorem{prop}{Proposition} [section]

\newtheorem{lemma}[prop]{Lemma}

\newtheorem{thm}[prop]{Theorem}

\newtheorem{cor}[prop]{Corollary}

\newtheorem*{thmA}{Proposition A}

\newtheorem*{thmB}{Proposition B}

\newtheorem*{thmC}{Proposition C}

\newtheorem*{thmD}{Theorem D}

\theoremstyle{definition}

\newtheorem{exa}{Example} [section]

\newtheorem*{claim}{Claim}

\newtheorem*{rmk}{Remark}

\newtheorem*{rmks}{Remarks}

\theoremstyle{remark}

\newtheorem*{pf}{Proof}

\newtheorem*{pfA}{Proof of Proposition A}

\newtheorem*{pfB}{Proof of Proposition B}

\newtheorem*{pfC}{Proof of Proposition C}

\newtheorem*{pfD}{Proof of Theorem D}

\newtheorem*{ack}{Acknowledgment}

\numberwithin{equation}{section}
\makeatletter
\@namedef{subjclassname@2020}{%
	\textup{2020} Mathematics Subject Classification}
\makeatother
\title{Flip Signatures}

\author{Sieye Ryu}

\email{Sieye Ryu <sieyeryu@gmail.com>}

\address{Institute of Mathematics and Statistics, University of Sao Paulo, (IME-USP), Rua do
Matao 1010, CEP 05508-090, Sao Paulo, Brazil}


\begin{abstract}
A $D_{\infty}$-topological Markov chain is a topological Markov chain provided with an action of the infinite dihedral group $D_{\infty}$. It is defined by two zero-one square matrices $A$ and $J$ satisfying $AJ=JA^{\textsf{T}}$ and $J^2=I$.
Flip signature is obtained from symmetric bilinear forms with respect to $J$ on the eventual kernel of $A$.
We modify Williams' decomposition theorem to prove flip signature is a $D_{\infty}$-conjugacy invariant.
We introduce natural $D_{\infty}$-actions on Ashley's eight-by-eight and the full two-shift. 
The Flip signatures show that Ashley's eight-by-eight and the full two-shift equipped with the natural $D_{\infty}$-actions are not $D_{\infty}$-conjugate.
We also discuss the notion of $D_{\infty}$-shift equivalence and the Lind zeta function. 
\end{abstract}

\maketitle

\section{Introduction}

A \textit{topological flip system} $(X, T, F)$ is a topological dynamical system $(X, T)$ consisting of a topological space $X$, a homeomorphism $T : X \rightarrow X$ and a topological conjugacy $F: (X, T^{-1}) \rightarrow (X, T)$ with $F^2 = \text{Id}_X$. (See the survey paper \cite{LR} for the further study of flip systems.)
We call the map $F$ a \textit{flip} for $(X, T)$.
If $F$ is a flip for a discrete-time topological dynamical system $(X, T)$, then the triple $(X, T, F)$ is called a \textit{$D_{\infty}$-system} because
the infinite dihedral group 
$$D_{\infty} = \langle a, b : ab=ba^{-1} \text{ and } b^2=1 \rangle$$
acts on $X$ as follows:
$$(a, x) \mapsto T(x) \qquad \text{and} \qquad (b, x) \mapsto F(x) \qquad (x \in X).$$

Two $D_{\infty}$-systems $(X, T, F)$ and  $(X', T', F')$ are said to be \textit{$D_{\infty}$-conjugate} if there is a $D_{\infty}$-equivariant homeomorphism $\theta: X \rightarrow X'$. In this case, we write
$$(X, T, F) \cong (X', T', F')$$ 
and call the map $\theta$ a \textit{$D_{\infty}$-conjugacy} from $(X, T, F)$ to $(X', T', F')$.

Suppose that $\mathcal{A}$ is a finite set. 
A \textit{topological Markov chain}, or TMC for short, $(\textsf{X}_A, \sigma_A)$ over $\mathcal{A}$ is a shift space which has a zero-one $\mathcal{A} \times \mathcal{A}$ matrix $A$ as a transition matrix:
$$\textsf{X}_A = \{x \in \mathcal{A}^{\mathbb{Z}} : A(x_i, x_{i+1}) = 1 \;\; \forall i \in \mathbb{Z} \}.$$ 
A $D_{\infty}$-system $(X, T, F)$ is said to be a \textit{$D_{\infty}$-}topological Markov chain, or \textit{$D_{\infty}$-}TMC for short,  if $(X, T)$ is a topological Markov chain.

Suppose that $(X, T)$ is a shift space.
A flip $F$ for $(X, T)$ is called a \textit{one-block flip} if $x_0 = x'_0$ implies $F(x)_0=F(x')_0$ for all $x$ and $x'$ in $X$. 
If $F$ is a one-block flip for $(X, T)$, then there is a unique map $\tau: \mathcal{A} \rightarrow \mathcal{A}$ such that
$$F(x)_i = \tau(x_{-i}) \qquad \text{and} \qquad \tau^2 = \text{Id}_{\mathcal{A}} \qquad (x \in X; i \in \mathbb{Z}).$$ 
Representation Theorem in \cite{KLP} says that if $(X, T, F)$ is a $D_{\infty}$-TMC, then there is a TMC $(X', T')$ and a one-block flip $F'$ for $(X', T')$ such that $(X, T, F)$ and $(X', T', F')$ are $D_{\infty}$-conjugate.

Suppose that $\mathcal{A}$ is a finite set and that $A$ and $J$ are zero-one $\mathcal{A} \times \mathcal{A}$ matrices satisfying
\begin{equation}\label{flip pair}
AJ=JA^{\textsf{T}} \qquad \text{and} \qquad J^2=I.
\end{equation}
Since $J$ is zero-one and $J^2=I$, it follows that $J$ is symmetric and that for any $a \in \mathcal{A}$, there is a unique $b \in \mathcal{A}$ such that $J(a, b)=1$. Thus, there is a unique permutation 
$\tau_J:\mathcal{A} \rightarrow \mathcal{A}$ of order two satisfying
$$J(a, b) =1 \qquad \Leftrightarrow \qquad \tau_J(a) = b \qquad (a, b \in \mathcal{A}).$$
It is obvious that the map $\varphi_{J} : \mathcal{A}^{\mathbb{Z}} \rightarrow \mathcal{A}^{\mathbb{Z}}$ defined by
$$ \varphi_J(x)_i = \tau_J(x_{-i}) \qquad (x \in X)$$ is a one-block flip for the full $\mathcal{A}$-shift $(\mathcal{A}^{\mathbb{Z}}, \sigma)$.
Since $AJ=JA^{\textsf{T}}$ implies
$$A(a, b) = A(\tau_J(b), \tau_J(a)) \qquad (a, b \in \mathcal{A}),$$
it follows that $\varphi_{J}(\textsf{X}_A) = \textsf{X}_A$. Thus, the restriction $\varphi_{A, J}$ of $\varphi_J$ to $\textsf{X}_A$ becomes a one-block flip for $(\textsf{X}_A, \sigma_A)$. 
A pair $(A, J)$ of zero-one $\mathcal{A} \times \mathcal{A}$ matrices satisfying
(\ref{flip pair})
will be called a \textit{flip pair}.

The classification of shifts of finite type up to conjugacy is one of the central problems in symbolic dynamics. There is an algorithm determining whether or not two one-sided shifts of finite type ($\mathbb{N}$-SFTs) are $\mathbb{N}$-conjugate. (See Section 2.1 in \cite{Ki}.) 
In the case of two-sided shifts of finite type
($\mathbb{Z}$-SFTs), however, one cannot determine whether or not two systems
are $\mathbb{Z}$-conjugate, even though many $\mathbb{Z}$-conjugacy invariants have been discovered. 
For instance, it is well known (Proposition 7.3.7 in \cite{LM}) that if two $\mathbb{Z}$-SFTs are $\mathbb{Z}$-conjugate, then their transition matrices have the same set of nonzero eigenvalues. In 1990, Ashley introduced an eight-by-eight zero-one matrix, which is called Ashley's eight-by-eight and asked whether or not it is $\mathbb{Z}$-conjugate to the full two-shift. (See Example 2.2.7 in \cite{Ki} or Section 3 in \cite{B}.) Since the characteristic polynomial of Ashley's eight-by-eight is $t^7(t-2)$, we could say Ashley's eight-by-eight is very simple in terms of spectrum and it is easy to prove that Ashley's eight-by-eight is not $\mathbb{N}$-conjugate to the full two-shift. Nevertheless, this problem has not been solved yet.
Meanwhile, both Ashley's eight-by-eight and the full-two shift have one-block flips. 
More precisely, if we set
$$A=\left[\begin{array}{rrrrrrrr} 1 & 1 & 0 & 0 & 0 & 0 & 0 & 0 \\ 0 & 0 & 1 & 0 & 0 & 0 & 1 & 0  \\ 0 & 0 & 0 & 1 & 0 & 1 & 0 & 0 \\ 0 & 1 & 0 & 0 & 0 & 0 & 0 & 1 \\ 1 & 0 & 0 & 0 & 1 & 0 & 0 & 0 \\ 0 & 0 & 0 & 0 & 1 & 0 & 0 & 1 \\ 0 & 0 & 1 & 0 & 0 & 1 & 0 & 0 \\ 0 & 0 & 0 & 1 & 0 & 0 & 1 & 0 \end{array}\right] \,\,  J=\left[\begin{array}{rrrrrrrr} 0 & 0 & 0 & 0 & 1 & 0 & 0 & 0 \\ 0 & 0 & 0 & 0 & 0 & 1 & 0 & 0  \\ 0 & 0 & 0 & 0 & 0 & 0 & 1 & 0 \\ 0 & 0 & 0 & 0 & 0 & 0 & 0 & 1 \\ 1 & 0 & 0 & 0 & 0 & 0 & 0 & 0 \\ 0 & 1 & 0 & 0 & 0 & 0 & 0 & 0 \\ 0 & 0 & 1 & 0 & 0 & 0 & 0 & 0 \\ 0 & 0 & 0 & 1 & 0 & 0 & 0 & 0 \end{array}\right],$$
\begin{equation}\label{matrices}
B = \left[\begin{array}{rr} 1 & 1 \\ 1 & 1 \end{array}\right], \;\; I = \left[\begin{array}{rr} 1 & 0 \\ 0 & 1 \end{array}\right] \;\; \text{and} \;\; K = \left[\begin{array}{rr} 0 & 1 \\ 1 & 0 \end{array}\right],
\end{equation}
then $A$ is Ashley's eight-by-eight, $\varphi_{A, J}$ is a unique one-block flip for $(\textsf{X}_A, \sigma_A)$, $B$ is the minimal zero-one matrix defining the full two-shift and $(\textsf{X}_B, \sigma_B)$ has exactly two one-block flips $\varphi_{B, I}$ and $\varphi_{B, K}$. 
It is natural to ask whether or not $(\textsf{X}_A, \sigma_A, \varphi_{A, J})$ is $D_{\infty}$-conjugate to $(\textsf{X}_B, \sigma_B, \varphi_{B, I})$ or $(\textsf{X}_B, \sigma_B, \varphi_{B, K})$. 
In this paper, we introduce the notion of \textit{flip signatures} and prove \begin{equation}\label{re1}
(\textsf{X}_A, \sigma_A, \varphi_{A, J}) \ncong (\textsf{X}_B, \sigma_B, \varphi_{B, I}),
\end{equation}
\begin{equation}\label{re2}
(\textsf{X}_A, \sigma_A, \varphi_{A, J}) \ncong (\textsf{X}_B, \sigma_B, \varphi_{B, K})\end{equation}
and
\begin{equation}\label{re3}
(\textsf{X}_B, \sigma_B, \varphi_{B, I}) \ncong (\textsf{X}_B, \sigma_B, \varphi_{B, K}).
\end{equation}

We first introduce analogues of elementary equivalence (EE), strong shift equivalence (SSE) and Williams' decomposition theorem for $D_{\infty}$-TMCs.
Let us recall the notions of EE and SSE. (See \cite{LM, W} for the details.) Suppose that $A$ and $B$ are zero-one square matrices. A pair $(D, E)$ of zero-one matrices satisfying
$$A=DE \qquad \text{and} \qquad B=ED$$ is said to be an \textit{elementary equivalence} (EE) \textit{from $A$ to $B$} and we write $(D, E): A \threesim B$. 
If $(D, E): A \threesim B$, then there is a $\mathbb{Z}$-conjugacy $\gamma_{D, E}$ from $(\textsf{X}_A, \sigma_A)$ to 
$(\textsf{X}_B, \sigma_B)$ satisfying
\begin{equation}\label{elementary conjugacy}
\gamma_{D, E}(x) = y \qquad \Leftrightarrow \qquad \forall \, i \in \mathbb{Z} \qquad D(x_i, y_i)=E(y_i, x_{i+1})=1.
\end{equation}
The map $\gamma_{D, E}$ is called an \textit{elementary conjugacy}. 

A \textit{strong shift equivalence} (SSE) \textit{of lag $l$ from $A$ to $B$} is a sequence of $l$ elementary equivalences 
$$(D_1, E_1): A \threesim A_1, \;\; (D_2, E_2): A_1 \threesim A_2, \;\; \cdots, \;\; (D_l, E_l): A_l \threesim B.$$
It is evident that if $A$ and $B$ are strong shift equivalent, then $(\textsf{X}_A, \sigma_A)$ and $(\textsf{X}_B, \sigma_B)$ are $\mathbb{Z}$-conjugate.
Williams' decomposition theorem, found in \cite{W}, says that 
every $\mathbb{Z}$-conjugacy between two $\mathbb{Z}$-TMCs can be decomposed into the composition of a finite number of elementary conjugacies.

To establish analogues of EE, SSE and Williams' decomposition theorem for $D_{\infty}$-TMCs, we first observe some properties of a $D_{\infty}$-system.
If $(X, T, F)$ is a $D_{\infty}$-system, 
then $(X, T, T^n \circ F)$ are also $D_{\infty}$-systems for all integers $n$. 
It is obvious that $T^n$ are $D_{\infty}$-conjugacies from $(X, T, F)$ to $(X, T, T^{2n} \circ F)$ and from $(X, T, T \circ F)$ to $(X, T, T^{2n+1} \circ F)$ for all integers $n$. 
For one's information, we will see that $(X, T, F)$ is not $D_{\infty}$-conjugate to $(X, T, T\circ F)$ in Proposition \ref{prop: 4.1}.

Let $(A, J)$ and $(B, K)$ be flip pairs. 
A pair $(D, E)$ of zero-one matrices satisfying
$$A=DE, \;\; B=ED \;\; \text{and} \;\; E=KD^{\textsf{T}}J$$ 
is said to be a \textit{$D_{\infty}$-half elementary equivalence} ($D_{\infty}$-HEE) from $(A, J)$ to $(B, K)$ and write $(D, E): (A, J) \threesim (B, K)$. 
In Proposition \ref{prop: 2.2}, we will see that
if $(D, E): (A, J) \threesim (B, K)$, then
the elementary conjugacy $\gamma_{D, E}$ from (\ref{elementary conjugacy}) becomes a $D_{\infty}$-conjugacy from $(\textsf{X}_A, \sigma_A, \varphi_{A, J})$ to $(\textsf{X}_B, \sigma_B, \sigma_B \circ \varphi_{B, K})$.
We call the map $\gamma_{D, E}$ a \textit{$D_{\infty}$-half elementary conjugacy from $(\textsf{X}_A$, $\sigma_A$, $\varphi_{A, J})$ to $(\textsf{X}_B$, $\sigma_B$, $\sigma_B \circ \varphi_{B, K})$}.

A sequence of $l$ $D_{\infty}$-half elementary equivalences
$$(D_1, E_1) : (A, J) \threesim (A_2, J_2), \;\; \cdots, \; \; (A_l, D_l) : (A_l, D_l) \threesim (B, K)$$
is said to be a \textit{$D_{\infty}$-strong shift equivalence} ($D_{\infty}$-SSE) of lag $l$ from $(A, J)$ to
$(B, K)$.
If there is a $D_{\infty}$-SSE of lag $l$ from $(A, J)$ to $(B, K)$, then 
$(\textsf{X}_A, \sigma_A, \varphi_{A, J})$ is $D_{\infty}$-conjugate to $(\textsf{X}_B, \sigma_B, \sigma_B^l \circ \varphi_{B, K})$.
If $l$ is an even number, then $(\textsf{X}_A, \sigma_A, \varphi_{A, J})$ is $D_{\infty}$-conjugate to $(\textsf{X}_B, \sigma_B, \varphi_{B, K})$, while if $l$ is an odd number, then $(\textsf{X}_A, \sigma_A, \varphi_{A, J})$ is $D_{\infty}$-conjugate to $(\textsf{X}_B, \sigma_B, \sigma_B \circ \varphi_{B, K})$.  
In Section 4, we will see that Williams' decomposition theorem can be modified as follows:

\begin{thmA}
Suppose that $(A, J)$ and $(B, K)$ are flip pairs.
\newline
\emph{(1)} Two $D_{\infty}$-TMCs $(\textsf{X}_A, \sigma_A, \varphi_{A, J})$ and $(\textsf{X}_B, \sigma_B, \varphi_{B, K})$ are $D_{\infty}$-conjugate if and only if there is a $D_{\infty}$-SSE of lag $2l$ between $(A, J)$ and $(B, K)$ for some positive integer $l$.
\newline
\emph{(2)} Two $D_{\infty}$-TMCs $(\textsf{X}_A, \sigma_A, \varphi_{A, J})$ and $(\textsf{X}_B, \sigma_B, \sigma_B \circ \varphi_{B, K})$ are $D_{\infty}$-conjugate if and only if there is a $D_{\infty}$-SSE of lag $2l-1$ between $(A, J)$ and $(B, K)$ for some positive integer $l$.
\end{thmA}

In order to introduce the notion of flip signatures, we discuss some properties of $D_{\infty}$-TMCs.
We first indicate notation.
If $\mathcal{A}_1$ and $\mathcal{A}_2$ are finite sets and $M$ is an $\mathcal{A}_1 \times \mathcal{A}_2$ zero-one matrix, then for each $a \in \mathcal{A}_1$,
we set
$$\mathcal{F}_M(a) = \{b \in \mathcal{A}_2: M(a, b) =1 \}$$
and for each $b \in \mathcal{A}_2$, we set
$$\mathcal{P}_M(b) = \{a \in \mathcal{A}_1: M(a, b) =1 \}.$$
When $(X, T)$ is a TMC,
we denote the set of all $n$-blocks occurring in points in $X$ by $\mathcal{B}_n(X)$ for all nonnegative integers $n$.

Suppose that $(A, J)$ and $(B, K)$ are flip pairs and that $(D, E)$ is a $D_{\infty}$-HEE from $(A, J)$ to $(B, K)$.
Since $B$ is zero-one and $B=ED$, it follows that
\begin{eqnarray}\label{eq1.3}
\mathcal{F}_D(a_1) \cap \mathcal{F}_D(a_2) \neq \varnothing \qquad &\Rightarrow & \qquad \mathcal{P}_E(a_1) \cap \mathcal{P}_E(a_2) = \varnothing \nonumber \\
\mathcal{P}_E(a_1) \cap \mathcal{P}_E(a_2) \neq \varnothing \qquad & \Rightarrow & \qquad
\mathcal{F}_D(a_1) \cap \mathcal{F}_D(a_2) = \varnothing,
\end{eqnarray}
for all $a_1, a_2 \in \mathcal{B}_1(\textsf{X}_A)$.
Suppose that $u$ and $v$ are real-valued functions defined on $\mathcal{B}_1(\textsf{X}_A)$ and  $\mathcal{B}_1(\textsf{X}_B)$, respectively.
If $|\mathcal{B}_1(\textsf{X}_A)|=m$ and $|\mathcal{B}_1(\textsf{X}_B)|=n$, then $u$ and $v$ can be regarded as vectors in $\mathbb{R}^m$ and $\mathbb{R}^n$, respectively. 
If $u$ and $v$ satisfy
\begin{equation}\label{key}
\forall \, a \in \mathcal{B}_1(\textsf{X}_A) \quad u(a) = \sum_{b \in \mathcal{F}_D(a)} v(b),
\end{equation}
then for each $a \in \mathcal{B}_1(\textsf{X}_A)$, we have
$$u(\tau_J(a)) u(a) = \sum_{b \in \mathcal{P}_E(a)} v(\tau_K(b)) \sum_{b \in \mathcal{F}_D(a)} v(b)$$
by $E = KD^{\textsf{T}}J$ and
(\ref{eq1.3}) leads to
$$\sum_{a \in \mathcal{B}_1(\textsf{X}_A)} u(\tau_J(a)) u(a) = \sum_{b \in \mathcal{B}_1(\textsf{X}_B)} \sum_{d \in \mathcal{P}_B(b)} v(\tau_K(d)) v(b).$$
Since $J$ and $K$ are symmetric, this formula can be expressed in terms of symmetric bilinear forms with respect to $J$ and $K$:
$$u^{\textsf{T}}Ju=(Bv)^{\textsf{T}}Kv.$$
We note that if both $A$ and $B$ have $\lambda$ as their real eigenvalues and $v$ is an eigenvector of $B$ corresponding to $\lambda$, then $u$ satisfying (\ref{key}) is an eigenvector of $A$ corresponding to $\lambda$.  
We consider the case where $A$ and $B$ have $0$ as their eigenvalues and find out some relationships between the symmetric bilinear forms of the generalized eigenvectors of $A$ and $B$ corresponding to $0$ when $(A, J)$ and $(B, K)$ are $D_{\infty}$-half elementary equivalent.

We call the subspace $\mathcal{K}(A)$ of $u \in \mathbb{R}^m$ such that $A^p u = 0$ for some $p \in \mathbb{N}$ the \textit{eventual kernel} of $A$:
$$\mathcal{K}(A) = \{u \in \mathbb{R}^m : A^p u = 0 \text{ for some } p \in \mathbb{N} \}.$$
If $u \in \mathcal{K}(A) \setminus \{0\}$ and $p$ is the smallest positive integer for which $A^pu=0$, then the ordered set 
$$\alpha = \{A^{p-1}u, \cdots, Au, u\}$$
is called a \textit{cycle of generalized eigenvectors of $A$ corresponding to $0$}. In this paper, we sometimes call $\alpha$ a \textit{cycle in $\mathcal{K}(A)$} for simplicity.
The vectors $A^{p-1}u$ and $u$ are called the \textit{initial vector} and the \textit{terminal vector} of $\alpha$, respectively and we write 
$$\text{Ini}(\alpha) = A^{p-1}u \qquad  \text{and} \qquad \text{Ter}(\alpha) = u.$$
We say that the length of $\alpha$ is $p$ and write $|\alpha|=p$.
It is well known \cite{FIS} that there is a basis for $\mathcal{K}(A)$ consisting of a union of disjoint cycles of generalized eigenvectors of $A$ corresponding to $0$.
The set of bases for $\mathcal{K}(A)$ consisting of a union of disjoint cycles of generalized eigenvectors of $A$ corresponding to $0$ is denoted by $\mathcal{B}as(\mathcal{K}(A))$.
We will prove the following proposition in Section 3.

\begin{thmB}
Suppose that $(D, E): (A, J) \threesim (B, K)$. 
Then there exist bases $\mathcal{E}(A) \in \mathcal{B}as(\mathcal{K}(A))$ and $\mathcal{E}(B) \in \mathcal{B}as(\mathcal{K}(B))$ such that
if $p>1$ and $\alpha = \{u_1, u_2, \cdots, u_p\}$ is a cycle in $\mathcal{E}(A)$ 
then one of the following holds.
\newline
\emph{(1)} There is a cycle $\beta=\{v_1, v_2, \cdots, v_{p+1}\}$ in $\mathcal{E}(B)$ 
such that
$$Dv_{k+1}=u_{k} \qquad \text{and} \qquad  Eu_k =v_k \qquad (k=1, \cdots, p).$$
\emph{(2)} There is a cycle $\beta = \{v_1, v_2, \cdots, v_{p-1}\}$ in $\mathcal{E}(B)$ 
such that
$$Dv_{k}=u_{k} \qquad \text{and} \qquad  Eu_{k+1} =v_k \qquad (k=1, \cdots, p-1).$$
In either case, we have
\begin{equation}\label{well sgn}
\emph{Ini}(\alpha)^{\textsf{T}}J\emph{Ter}(\alpha) = \emph{Ini}(\beta)^{\textsf{T}}K\emph{Ter}(\beta).
\end{equation}
\end{thmB}

In Lemma \ref{generalized eigenvectors}, we will show that there is a basis $\mathcal{E}(A) \in \mathcal{B}as(\mathcal{K}(A))$ such that
the left-hand side of (\ref{well sgn}) is not $0$ for every cycle $\alpha$ in $\mathcal{E}(A)$.
In this case, we define the sign of a cycle $\alpha = \{u_1, u_2, \cdots, u_p\}$ in $\mathcal{E}(A)$ by
$$\text{sgn}(\alpha) = \begin{cases}+1 \qquad \text{if } u_1^{\textsf{T}}Ju_p>0 \\ -1 \qquad \text{if } u_1^{\textsf{T}}Ju_p<0. \end{cases}$$
We denote the set of $|\alpha|$ such that $\alpha$ is a cycle in $\mathcal{E}(A)$ by $\mathcal{I}nd(\mathcal{K}(A))$.
It is clear that $\mathcal{I}nd(\mathcal{K}(A))$ is independent of the choice of basis for $\mathcal{K}(A)$.
We denote the union of the cycles $\alpha$ of length $p$ in $\mathcal{E}(A)$ by
$\mathcal{E}_p(A)$  for each $p\in \mathcal{I}nd(\mathcal{K}(A))$
and define the sign of $\mathcal{E}_p(A)$ by
$$\text{sgn}(\mathcal{E}_p(A)) = \prod_{\{\alpha : \alpha \text{ is a cycle in } \mathcal{E}_p(A)\}} \text{sgn}(\alpha).$$
In Section 3, we will prove the sign of $\mathcal{E}_p(A)$ is also independent of the choice of basis for $\mathcal{K}(A)$ if it is well-defined.

\begin{thmC}
Suppose that $\mathcal{E}(A)$ and $\mathcal{E}'(A)$ are two distinct bases in $\mathcal{B}as(\mathcal{K}(A))$
such that the sign of every cycle in both $\mathcal{E}(A)$ and $\mathcal{E}'(A)$ is well-defined.
For each $p \in \mathcal{I}nd(\mathcal{K}(A))$, we have
$$\emph{sgn}(\mathcal{E}_p(A))=\emph{sgn}(\mathcal{E}'_p(A)).$$
\end{thmC}

Suppose that $\mathcal{E}(A) \in \mathcal{B}as(\mathcal{K}(A))$ 
and that the sign of every cycle in $\mathcal{E}(A)$ is well-defined.
We arrange the elements of $\mathcal{I}nd(\mathcal{K}(A)) = \{p_1, p_2, \cdots, p_A\}$  
to satisfy 
$$p_1 < p_2 < \cdots p_A$$
and write $$\varepsilon_{p}=\text{sgn}(\mathcal{E}_p(A)).$$
If $|\mathcal{I}nd(\mathcal{K}(A))|=k$, then
the $k$-tuple $(\varepsilon_{p_1}, \varepsilon_{p_2}, \cdots, \varepsilon_{p_A})$ is called the \textit{flip signature of $(A, J)$} and $\varepsilon_{p_A}$ is called the \textit{leading signature of $(A, J)$}.
The flip signature of $(A, J)$ is denoted by 
$$\emph{F.Sig}(A, J) = (\varepsilon_{p_1}, \varepsilon_{p_2}, \cdots, \varepsilon_{p_A}).$$
The following is the main result of this paper.

\begin{thmD}
Suppose that $(A, J)$ and $(B, K)$ are flip pairs and that $(\textsf{X}_A, \sigma_A, \varphi_{A, J})$ and $(\textsf{X}_B, \sigma_B, \varphi_{B, K})$ are $D_{\infty}$-conjugate. 
If 
$$\text{F.Sig}(A, J) = (\varepsilon_{p_1}, \varepsilon_{p_2}, \cdots, \varepsilon_{p_A})$$
and
$$\text{F.Sig}(B, K) = (\varepsilon_{q_1}, \varepsilon_{q_2}, \cdots, \varepsilon_{q_B}),$$
then $\text{F.Sig}(A, J)$ and $\text{F.Sig}(B, K)$ have the same number of $-1$'s and
$$\epsilon_{p_A} = \epsilon_{q_B}.$$
\end{thmD}

In Section 7, we will compute the flip signatures of $(A, J)$, $(B, I)$ and $(B, K)$, where $A$, $J$, $B$, $I$ and $K$ are as in (\ref{matrices}) and prove (\ref{re1}), (\ref{re2}) and (\ref{re3}).
Actually, we can obtain (\ref{re1}), (\ref{re2}) and (\ref{re3}) from the Lind zeta functions.
In \cite{KLP}, an explicit formula for the Lind zeta function for a $D_{\infty}$-TMC was established, which can be expressed in terms of matrices from flip pairs. From its formula (See also Section 6.), it is obvious that the Lind zeta function is a $D_{\infty}$-conjugacy invariant. 
In Example 7.1, we will see that the Lind zeta functions of $(\textsf{X}_A, \sigma_A, \varphi_{A, J})$, $(\textsf{X}_B, \sigma_B, \varphi_{B, I})$ and $(\textsf{X}_B, \sigma_B, \varphi_{B, K})$ are all different.  
In Section 6, we introduce the notion of $D_{\infty}$-shift equivalence ($D_{\infty}$-SE) which is an analogue of shift equivalence and prove that $D_{\infty}$-SE is a $D_{\infty}$-conjugacy invariant. In Example 7.2, we will see that there are $D_{\infty}$-SEs between $(A, J)$, $(B, I)$ and $(B, K)$ pairwise. 
So the existence of $D_{\infty}$-shift equivalence between two flip pairs does not imply that the corresponding $\mathbb{Z}$-TMCs share the same Lind zeta functions. This is a contrast to the fact that the existence of shift equivalence between two defining matrices $A$ and $B$ implies the coincidence of the Artin-Mazur zeta functions of the $\mathbb{Z}$-TMCs $(\textsf{X}_A, \sigma_A)$ and $(\textsf{X}_B, \sigma_B)$.
Example 7.4, however, says that
the coincidence of the Lind zeta functions of two $D_{\infty}$-TMCs does not guarantee the existence of $D_{\infty}$-shift equivalence between their flip pairs. This is analogous to the case of $\mathbb{Z}$-TMCs because the coincidence of the Artin-Mazur zeta functions of two $\mathbb{Z}$-TMCs does not guarantee the existence of SE between their defining matrices. (See Section 7.4 in \cite{LM}.)

This paper is organized as follows.
In Section 2, we introduce the notions of $D_{\infty}$-half elementary equivalence and $D_{\infty}$-strong shift equivalence.
In Section 3, we investigate symmetric bilinear forms with respect to $J$ and $K$ on the eventual kernels of $A$ and $B$ when two flip pairs $(A, J)$ and $(B, K)$ are $D_{\infty}$-half elementary equivalent. In the same section, we prove Proposition B and Proposition C. 
Proposition A and Theorem D will be proved in Section 4 and Section 5, respectively.
In section 6, we discuss the notion of $D_{\infty}$-shift equivalence and the Lind zeta function.
Section 7 consists of examples.

\begin{ack}
The author gratefully acknowledge the support of FAPESP (Grant No. 2018/12482-3) and Institute of Mathematics and Statistics in University of Sao Paulo (IME-USP).
\end{ack}

\section{$D_{\infty}$-Strong Shift Equivalence}

\label{sec:second}

Let $(A, J)$ and $(B, K)$ be flip pairs. A pair $(D, E)$ of zero-one matrices satisfying
$$A=DE, \quad B=ED, \quad \text{and} \quad E=KD^{\textsf{T}}J$$
is said to be a \textit{$D_{\infty}$-half elementary equivalence from $(A, J)$ to $(B, K)$}. 
If there is a $D_{\infty}$-half elementary equivalence from $(A, J)$ to $(B, K)$, then we write $(D, E) : (A, J) \threesim (B, K)$.
We note that symmetricities of $J$ and $K$ imply 
$$E=KD^{\textsf{T}}J \qquad \Leftrightarrow \qquad D=JE^{\textsf{T}}K.$$ 

\begin{prop}
\label{prop: 2.2}
If $(D, E) : (A, J) \threesim (B, K)$, then $(\textsf{X}_A, \sigma_A, \varphi_{J, A})$ is $D_{\infty}$-conjugate to $(\textsf{X}_B, \sigma_B, \sigma_B \circ \varphi_{K, B})$.
\end{prop}

\begin{pf}
Since $D$ and $E$ are zero-one and $A=DE$, it follows that for all $a_1a_2 \in \mathcal{B}_2(\textsf{X}_A)$, there is a unique $b\in \mathcal{B}_1(\textsf{X}_B)$ such that 
$$D(a_1, b) = E(b, a_2)=1.$$
We denote the block map which sends $a_1a_2 \in \mathcal{B}_2(\textsf{X}_A)$ to $b \in \mathcal{B}_1(\textsf{X}_B)$ by $\Gamma_{D, E}$. 
If we define the map $\gamma_{D, E} : (\textsf{X}_A, \sigma_A) \rightarrow (\textsf{X}_B, \sigma_B)$ by
$$\gamma_{D, E}(x)_i = \Gamma_{D, E} \left(x_i x_{i+1} \right) \qquad (x \in \textsf{X}_A; \; i \in \mathbb{Z}),$$
then we have $\gamma_{D, E} \circ \sigma_A = \sigma_B \circ \gamma_{D, E}$.

Since $(E, D): (B, K) \threesim (A, J)$, we can define the block map $\Gamma_{E, D} : \mathcal{B}_2(\textsf{X}_B) \rightarrow \mathcal{B}_1(\textsf{X}_A)$ and the map $\gamma_{E, D}: (\textsf{X}_B, \sigma_B) \rightarrow (\textsf{X}_A, \sigma_A)$ in the same way.
Since $\gamma_{E, D} \circ \gamma_{D, E} = \text{Id}_{\textsf{X}_A}$ and $\gamma_{D, E} \circ \gamma_{E, D} = \text{Id}_{\textsf{X}_B}$, it follows that $\gamma_{D, E}$ is one-to-one and onto.

It remains to show that 
\begin{equation}\label{eq: 2.5}
\gamma_{D, E} \circ \varphi_{A, J} = \left( \sigma_B \circ \varphi_{B, K} \right) \circ \gamma_{D, E}.
\end{equation}
Since $E = KD^{\textsf{T}}J$, it follows that 
$$E(b, a) = 1 \quad \Leftrightarrow \quad D(\tau_J(a), \tau_K(b))=1 \qquad (a \in \mathcal{B}_1(\textsf{X}_A),  b\in \mathcal{B}_1(\textsf{X}_B)).$$
This is equivalent to
$$D(a, b) = 1 \quad \Leftrightarrow \quad E(\tau_K(b), \tau_J(a))=1 \qquad (a \in \mathcal{B}_1(\textsf{X}_A),  b\in \mathcal{B}_1(\textsf{X}_B)).$$
Thus, we obtain
\begin{equation}\label{eq: 2.6}
\Gamma_{D, E}(a_1 a_2) =b \quad \Leftrightarrow \quad \Gamma_{D, E} \left( \tau_J(a_2) \tau_J(a_1) \right) = \tau_K(b) \qquad \left( a_1 a_2 \in \mathcal{B}_2 (\textsf{X}_A) \right).
\end{equation}
By (\ref{eq: 2.6}), we have
\begin{eqnarray*}
\gamma_{D, E} \circ \varphi_{J, A} (x)_i
&=& \Gamma_{D, E}(\tau_J(x_{-i}) \tau_J(x_{-i-1}))
= \tau_K (\Gamma_{D, E}(x_{-i-1}x_{-i}))\\
&=& \varphi_{B, K} \circ \gamma_{D, E} (x)_{i+1}
= \left( \sigma_B \circ \varphi_{B, K} \right) \circ \gamma_{D, E} (x)_i 
\end{eqnarray*}
for all $x \in \textsf{X}_A$ and $i \in \mathbb{Z}$ and this proves (\ref{eq: 2.5}). 
\hfill $\Box$
\end{pf}

Let $(A, J)$ and $(B, K)$ be flip pairs. 
A sequence of $l$ half elementary equivalences 
$$(D_1, E_1) : (A, J) \threesim (A_2, J_2),$$
$$(D_2, E_2) : (A_2, J_2) \threesim (A_3, J_3),$$
$$\vdots$$
$$(D_l, E_l): (A_l, J_l) \threesim (B, K)$$ 
is said to be a  \textit{$D_{\infty}$-SSE of lag $l$ from $(A, J)$ to $(B, K)$}.
If there is a $D_{\infty}$-SSE of lag $l$ from $(A, J)$ to $(B, K)$, then we say that $(A, J)$ is $D_{\infty}$-strong shift equivalent to $(B, K)$ and write $(A, J) \approx (B, K)$ (lag $l$).

By Proposition \ref{prop: 2.2}, we have
\begin{equation}\label{eq: 2.7}
(A, J) \approx (B, K) \; (\text{lag} \; l) \quad \Rightarrow \quad (\textsf{X}_A, \sigma_A, \varphi_{J, A}) \cong (\textsf{X}_B, \sigma_B, {\sigma_B}^l \circ \varphi_{K, B}).
\end{equation}
Because ${\sigma_B}^l$ is a conjugacy from $(\textsf{X}_B, \sigma_B, \varphi_{K, B})$ to $(\textsf{X}_B, \sigma_B, {\sigma_B}^{2l} \circ \varphi_{K, B})$, the implication in (\ref{eq: 2.7}) can be rewritten as follows:
\begin{equation}\label{eq: 2.8}
(A, J) \approx (B, K) \; (\text{lag} \; 2l) \quad \Rightarrow \quad (\textsf{X}_A, \sigma_A, \varphi_{J, A}) \cong (\textsf{X}_B, \sigma_B, \varphi_{K, B})
\end{equation}
and 
\begin{equation}\label{eq: 2.9}
(A, J) \approx (B, K) \; (\text{lag} \; 2l-1) \quad \Rightarrow \quad (\textsf{X}_A, \sigma_A, \varphi_{J, A}) \cong (\textsf{X}_B, \sigma_B, \sigma_B \circ \varphi_{K, B}).
\end{equation}  
In Section \ref{proof of thm1}, we will prove Proposition A which says that the converses of (\ref{eq: 2.8}) and (\ref{eq: 2.9}) are also true.

\section{Symmetric Bilinear Forms}
Suppose that $(A, J)$ is a flip pair and that $|\mathcal{B}_1(\textsf{X}_A)|=m$.
Let $V$ be an $m$-dimensional vector space over the filed $\mathbb{C}$ of complex numbers.
Let $\langle u, v \rangle_J$ denote the bilinear form $V \times V \rightarrow \mathbb{C}$ defined by 
$$(u, v)\mapsto u^{\textsf{T}}J\bar{v} \qquad (u, v \in V).$$ 
Since $J$ is a non-singular symmetric matrix, it follows that the bilinear form $\langle \;\;\; , \;\;\; \rangle_J$ is symmetric and non-degenerate.
If $u, v \in V$ and $\langle u, v \rangle_J = 0$, then $u$ and $v$ are said to be \textit{orthogonal with respect to $J$} and we write $u \perp_J v$.
From $AJ=JA^{\textsf{T}}$, we see that $A$ itself is the adjoint of $A$ in the following sense:
\begin{equation}\label{6}
\langle A{u, v} \rangle_J = \langle u, Av \rangle_J.
\end{equation}
If $\lambda$ is an eigenvalue of $A$ and $u$ is an eigenvector of $A$ corresponding to $\lambda$, then for any $v \in V$, we have
\begin{equation}\label{7}
\lambda \langle {u}, {v} \rangle_J = \langle  \lambda{u}, {v} \rangle_J = \langle  A{u}, {v} \rangle_J = \langle {u}, A{v} \rangle_J.
\end{equation}

Let $\text{sp}(A)$ denote the set of eigenvalues of $A$.
For each $\lambda \in \text{sp}(A)$, let $\mathcal{K}_{\lambda}(A)$ denote the set of $u \in V$ such that $(A-\lambda I)^p u =0$ for some $p \in \mathbb{N}$:
$$\mathcal{K}_{\lambda}(A) = \{u \in V : \exists \, p \in \mathbb{N} \text{ s.t. } (A-\lambda I)^p {u} = 0 \}.$$
If $u \in \mathcal{K}_{\lambda}(A) \setminus \{0\}$ and $p$ is the smallest positive integer for which $(A-\lambda I)^pu=0$, then the ordered set 
$$\alpha = \{(A-\lambda I)^{p-1}u, \cdots, (A-\lambda I)u, u\}$$
is called a \textit{cycle of generalized eigenvectors of $A$ corresponding to $\lambda$}. 
The vectors $(A-\lambda I)^{p-1}u$ and $u$ are called the \textit{initial vector} and the \textit{terminal vector} of $\alpha$, respectively and we write 
$$\text{Ini}(\alpha) = (A-\lambda I)^{p-1}u \qquad  \text{and} \qquad \text{Ter}(\alpha) = u.$$
We say that the length of $\alpha$ is $p$ and write $|\alpha|=p$.
It is well known \cite{FIS} that there is a basis for $\mathcal{K}_{\lambda}(A)$ consisting of a union of disjoint cycles of generalized eigenvectors of $A$ corresponding to $\lambda$.
From here on, when we say $\alpha=\{u_1, \cdots, u_p\}$ is a cycle in $\mathcal{K}_{\lambda}(A)$, it means $\alpha$ is a cycle of generalized eigenvectors of $A$ corresponding to $\lambda$, $\text{Ini}(\alpha) = u_1$, $\text{Ter}(\alpha)=u_p$ and $|\alpha|=p$.

Suppose that $\mathcal{U}(A)$ is a basis for $V$ consisting of generalized eigenvectors of $A$ and that $\mathcal{E}(A)$ is the subset of $\mathcal{U}(A)$ consisting of the generalized eigenvectors of $A$ corresponding to $0$.
Non-degeneracy of $\langle \;\;\;  , \;\;\; \rangle_J$ says that for each $u \in \mathcal{E}(A)$, there is a $v \in \mathcal{U}(A)$ such that $\langle u, v \rangle_J \neq 0$. 
The following lemma says that the vector $v$ must be in $\mathcal{E}(A)$.

\begin{lemma}\label{distinct eigenvalues}
Suppose that $\lambda, \mu \in \emph{sp}(A)$.
If $\lambda$ is distinct from the complex conjugate $\bar{\mu}$ of $\mu$, then
$\mathcal{K}_{\lambda}(A) \perp_J \mathcal{K}_{\mu}(A)$.
\end{lemma}

\begin{pf}
Suppose that 
$$\alpha = \{u_1, \cdots, u_p\} \qquad \text{and} \qquad \beta = \{v_1, \cdots, v_q\}$$
are cycles in $\mathcal{K}_{\lambda}(A)$ and $\mathcal{K}_{\mu}(A)$, respectively.
Since (\ref{7}) implies
$$\lambda \langle u_1, v_1 \rangle_J = \langle u_1, Av_1 \rangle_J = \bar{\mu} \langle u_1, v_1 \rangle_J, $$
it follows that 
$$\langle u_1, v_1 \rangle_J = 0.$$ 
Using (\ref{7}) again, we have
$$\lambda \langle u_1, v_{j+1} \rangle_J =  \langle u_1, \mu v_{j+1} + v_{j} \rangle_J  = \bar{\mu} \langle u_1, v_{j+1} \rangle_J + \langle u_1, v_j \rangle_J$$
for each $j=1, \cdots, q-1$.
By mathematical induction on $j$, we see that 
$$\langle u_1, v_j \rangle_J = 0 \qquad (j=1, \cdots, q).$$
Applying the same process to each $u_2, \cdots u_p$, we obtain
$$\forall \, i=1, \cdots, p, \quad \forall \, j=1, \cdots, q \qquad \langle u_i, v_j \rangle_J = 0.$$
\hfill $\Box$
\end{pf}

\begin{rmk} Suppose that $\mathcal{E}(A)=\{u_1, \cdots, u_p\}$ is a basis for $\mathcal{K}_0(A)$ consisting of generalized eigenvectors of $A$ corresponding to $0$.
If $T$ is the $m \times p$ matrix whose $i$-th column is $u_i$ for each $i=1, \cdots, p$, then non-degeneracy of $\langle \;\;, \;\; \rangle_J$ and Lemma \ref{distinct eigenvalues} implies that
$T^{\textsf{T}}JT$ is non-singular.	
\end{rmk}

From here on, we restrict our attention to the zero eigenvalue and the generalized eigenvectors corresponding to $0$. For notational simplicity, the smallest subspace of $V$ containing all generalized eigenvectors of $A$ corresponding to $0$ is denoted by $\mathcal{K}(A)$ and we call the subspace
$\mathcal{K}(A)$ of $V$ the \textit{eventual kernel} of $A$. 
We may assume that the eventual kernel of $A$ is a real vector space.
The set of bases for $\mathcal{K}(A)$ consisting of a union of disjoint cycles of generalized eigenvectors of $A$ corresponding to $0$ is denoted by $\mathcal{B}as(\mathcal{K}(A))$.
If $\mathcal{E}(A) \in \mathcal{B}as(\mathcal{K}(A))$, the set of $|\alpha|$ such that 
$\alpha$ is a cycle in $\mathcal{E}(A)$ is denoted by $\mathcal{I}nd(\mathcal{K}(A))$
and we call $\mathcal{I}nd(\mathcal{K}(A))$ the \textit{index set for the eventual kernel of $A$}.
It is clear that $\mathcal{I}nd(\mathcal{K}(A))$ is independent of the choice of $\mathcal{E}(A) \in \mathcal{B}as(\mathcal{K}(A))$.
When $p \in \mathcal{I}nd(\mathcal{K}(A))$,
we denote the union of the cycles of length $p$ in $\mathcal{E}(A)$ by $\mathcal{E}_p(A)$.

\begin{lemma} \label{non-singular}
Suppose that $\mathcal{E}(A) \in \mathcal{B}as(\mathcal{K}(A))$ and that $p \in  \mathcal{I}nd(\mathcal{K}(A))$.
If
$\mathcal{E}_p(A)=\{u_1, \cdots, u_r\}$
and $T_p$ is the $m \times r$ matrix whose $i$-th column is $u_i$ for each $i=1, \cdots, r$, then $T_p^{\textsf{T}}JT_p$ is non-singular.
\end{lemma}
\begin{pf}	
We only consider the case where $ \mathcal{I}nd(\mathcal{K}(A))=\{p, q\}$ $(p<q)$ and
both $\mathcal{E}_p(A)$ and $\mathcal{E}_q(A)$ have one cycles. 
If $T$ is the $m \times (p+q)$ matrix defined by
$$T =\left[\begin{array}{cc} T_p & T_q  \end{array}\right],$$
then 
$$T^{\textsf{T}}JT = \left[\begin{array}{cc} T_p^{\textsf{T}}JT_p &  T_p^{\textsf{T}}JT_q \\ T_q^{\textsf{T}}JT_p &  T_q^{\textsf{T}}JT_q \end{array}\right]$$ 
is non-singular by remark of Lemma \ref{distinct eigenvalues}.

Suppose that $\alpha=\{u_1, \cdots, u_p\}$ and $\beta = \{v_1, \cdots,v_q\}$ are cycles in $\mathcal{E}_p(A)$ and $\mathcal{E}_q(A)$, respectively. 
By (\ref{6}), we have 
$$ \langle u_1, u_{i} \rangle_J = \langle u_1, Au_{i+1} \rangle_J = \langle Au_1, u_{i+1} \rangle_J = 0$$
and 
\begin{equation}\label{computation}
\langle u_{i+1}, u_{j} \rangle_J - \langle u_{i}, u_{j+1} \rangle_J = \langle u_{i+1}, Au_{j+1} \rangle_J - \langle u_{i}, u_{j+1} \rangle_J=0
\end{equation}
for each $i, j=1, \cdots, p-1$.
Thus, if we set $\langle u_i, u_{p} \rangle_J = b_i$ for each $i=1, 2, \cdots, p$, then
$T_p^{\textsf{T}}JT_p$ is of the form
$$T_p^{\textsf{T}}JT_p = \left[\begin{array}{ccccccc} 0 & 0 & 0 & \cdots & 0 & 0 & b_1 \\ 0 & 0 & 0 & \cdots & 0 & b_1 & b_2  \\ 0 & 0 & 0 & \cdots & b_1 & b_2 & b_3 \\ \vdots & \vdots & \vdots &  & \vdots & \vdots & \vdots \\ b_1 & b_2 & b_3 & \cdots & b_{p-2} & b_{p-1} & b_p  \end{array}\right].$$
Obviously, $T_q^{\textsf{T}}JT_q$ is of the same form. We set
$$T_q^{\textsf{T}}JT_q = \left[\begin{array}{ccccccc} 0 & 0 & 0 & \cdots & 0 & 0 & d_1 \\ 0 & 0 & 0 & \cdots & 0 & d_1 & d_2  \\ 0 & 0 & 0 & \cdots & d_1 & d_2 & d_3 \\ \vdots & \vdots & \vdots &  & \vdots & \vdots & \vdots \\ d_1 & d_2 & d_3 & \cdots & d_{q-2} & d_{q-1} & d_q  \end{array}\right].$$
Now we consider $T_p^{\textsf{T}}JT_q$. 
By (\ref{6}) again, we have
\begin{eqnarray*}
\langle u_1, v_k \rangle_J &=& 0 \qquad (k=1, \cdots, q-1),\\
\langle u_2, v_k \rangle_J &=& 0 \qquad (k=1, \cdots, q-2),\\
\vdots&& \\
\langle u_p, v_k \rangle_J &=& 0 \qquad (k=1, \cdots, q-p).
\end{eqnarray*}
If we set $\langle u_i, v_q \rangle_J = c_{i}$ for each $i= 1, 2, \cdots, p$, then the argument in (\ref{computation}) shows that
$T_p^{\textsf{T}}JT_q$ is of the form
$$T_p^{\textsf{T}}JT_q = \left[\begin{array}{ccccccccc} 0 & \cdots & 0 & 0 & 0 & \cdots & 0 & 0 & c_1 \\ 0 & \cdots & 0 & 0 & 0 & \cdots & 0 & c_1 & c_2  \\ 0 & \cdots & 0 & 0 & 0 & \cdots & c_1 & c_2 & c_3 \\ \vdots & & \vdots & \vdots & \vdots &  & \vdots & \vdots & \vdots \\ 0 & \cdots & c_1 & c_2 & c_3 & \cdots & c_{p-2} & c_{p-1} & c_p  \end{array}\right].$$
Finally, $T_q^{\textsf{T}}JT_p$ is the transpose of $T_p^{\textsf{T}}JT_q$.
Hence, 
$b_1$ and $d_1$ must be nonzero and
we have
$\text{Rank}(T_p^{\textsf{T}}JT_p)=p$
and 
$\text{Rank}(T_q^{\textsf{T}}JT_q)=q$.
\hfill $\Box$
\end{pf}

The aim of this section is to find out a relationship between $\langle \;\; , \;\; \rangle_J$ and $\langle \;\; , \;\; \rangle_K$ on bases $\mathcal{E}(A) \in \mathcal{B}as(\mathcal{K}(A))$ and $\mathcal{E}(B) \in \mathcal{B}as(\mathcal{K}(B))$ when $(D, E): (A, J) \threesim (B, K)$. The following lemma will provide us good bases to handle.

\begin{lemma}\label{generalized eigenvectors}
Suppose that $A$ has the zero eigenvalue.
There is a basis $\mathcal{E}(A) \in \mathcal{B}as(\mathcal{K}(A))$ having the following properties.
\newline
\emph{(1)} If $\alpha$ is a cycle in $\mathcal{E}(A)$, then
$$\langle \emph{Ini}(\alpha), \emph{Ter}(\alpha) \rangle_J \neq 0.$$
\emph{(2)} Suppose that $\alpha$ is a cycle in $\mathcal{E}(A)$ with $\text{Ter}(\alpha)=u$ and $|\alpha|=p$. For each $k = 0, 1, \cdots p-1$, $v=A^{p-1-k}u$ is the unique vector in $\alpha$ such that $\langle A^ku, v \rangle_J \neq 0$.
\newline
\emph{(3)} If $\alpha$ and $\beta$ are distinct cycles in $\mathcal{E}(A)$, then
$$\emph{span}(\alpha) \perp_J \emph{span}(\beta).$$  
\end{lemma}

\begin{pf}
(1) First we consider the case where $\mathcal{E}(A)$ has only one cycle $\alpha=\{u_1, \cdots, u_p \}$.
By (\ref{6}), we have 
\begin{equation}\label{orthogonal}
\langle u_1, u_{i} \rangle_J = \langle u_1, Au_{i+1} \rangle_J = \langle Au_1, u_{i+1} \rangle_J = 0 \qquad (i=1, \cdots, p-1).
\end{equation}
By non-degeneracy of  $\langle \;\;\; , \;\;\; \rangle_J$ and Lemma \ref{distinct eigenvalues}, 
$\langle {u}_1, {u}_p \rangle_J $ must be nonzero. 

Suppose that $\mathcal{E}(A)$ is the union of disjoint cycles $\alpha_1, \cdots, \alpha_r$ of generalized eigenvectors of $A$ corresponding to $0$ for some $r>1$ and that 
$|\alpha_1| \leq |\alpha_2| \leq \cdots \leq |\alpha_r|$.
Assuming  
$$\langle \text{Ini}(\alpha_j), \text{Ter}(\alpha_j) \rangle_J \neq 0 \qquad (j=1, \cdots, r-1),$$ 
we will construct a cycle $\beta$ of generalized eigenvectors of $A$ corresponding to $0$ such that
the union of the cycles $\alpha_1, \cdots, \alpha_{r-1}, \beta$ forms a basis for $\mathcal{K}(A)$ and that $\langle \text{Ini}(\beta), \text{Ter}(\beta) \rangle_J \neq 0.$

If we set $\alpha_r = \{w_1, \cdots, w_q\}$,
the argument used in (\ref{orthogonal}) shows that
\begin{equation} \label{eq1234}
\langle w_1, w_{j} \rangle_J = 0 \qquad (j=1, \cdots, q-1).
\end{equation}
Suppose that $\alpha = \{u_1, \cdots, u_p \}$
is a cycle in $\mathcal{E}(A)$ which is distinct from $\alpha_r$. 
If $|\alpha|=p<q$,
then we have 
$$\langle w_1, A^ju_p \rangle_J = \langle A^{q-1}w_q , A^ju_p \rangle_J = \langle w_q, A^{j+q-1}u_p \rangle_J = 0$$ 
for all $j=0, \cdots, p-1$.
From non-degeneracy of $
\langle \;\; , \;\; \rangle_J$, Lemma \ref{distinct eigenvalues} and (\ref{eq1234}), 
it follows that 
$$|\alpha_1| \leq |\alpha_2| \leq \cdots \leq |\alpha_{r-1}| < |\alpha_r| \qquad \Rightarrow \qquad \langle w_1, w_q \rangle_J \neq 0.$$
In the case where there are other cycles of length $q$ in $\mathcal{E}(A)$, however, $\langle w_1, w_q \rangle_J$ is not necessarily nonzero. 
When $\langle w_1, w_q \rangle_J = 0$, there is a vector $v \in \mathcal{E}(A)$ such that $\langle w_1, v \rangle_J \neq 0$ by non-degeneracy of  $\langle \;\;\; , \;\;\; \rangle_J$ and Lemma \ref{distinct eigenvalues}. 
Since $\langle w_1, v \rangle_J = \langle w_q, A^{q-1}v \rangle_J$, 
it follows that
$v$ must be the terminal vector of a cycle in $\mathcal{E}(A)$ of length $q$ by the maximality of $q$. 
We put $v_1 = A^{q-1}v$ and $v_q=v$ 
and find a number $k \in \mathbb{R}\setminus \{0\}$ such that $\langle {w}_1-kv_1, w_q-kv_q \rangle_J \neq 0$.
We denote the cycle whose terminal vector is $w_q-kv_q$ by $\beta$. It is obvious that the length of $\beta$ is $q$ and that the union of the cycles $\alpha_1, \cdots, \alpha_{r-1}, \beta$ forms a basis of $\mathcal{K}(A)$.

(2) We assume that $\mathcal{E}(A)$ has property (1) and that $\alpha=\{u_1, \cdots, u_p\}$ is a cycle in $\mathcal{E}(A)$.
By (\ref{orthogonal}), we have $\langle u_1, u_i \rangle_J =0$ for all $i=1, \cdots, p-1$.
By (\ref{6}), we have
\begin{equation}\label{zerovalue}
\langle u_{i+1}, u_{j} \rangle_J - \langle u_{i}, u_{j+1} \rangle_J = \langle u_{i+1}, Au_{j+1} \rangle_J - \langle u_{i}, u_{j+1} \rangle_J=0
\end{equation}
for each $i, j = 1, \cdots, p-1$.
%and this implies 
%\begin{eqnarray}\label{zerovalue}
%\langle u_{i+1}, u_{j} \rangle_J = \langle u_{i}, u_{j+1} \rangle_J  \qquad (i , j = 1, 2 \cdots, p-1).
%\end{eqnarray}
Let $\langle u_i, u_{p} \rangle_J = b_i$ for each $i=1, 2, \cdots, p$.
If $T_{\alpha}$ is the $m \times p$ matrix whose $i$-th column is $u_i$, then $T_{\alpha}^{\textsf{T}}JT_{\alpha}$ is of the form
$$T_{\alpha}^{\textsf{T}}JT_{\alpha} = \left[\begin{array}{ccccccc} 0 & 0 & 0 & \cdots & 0 & 0 & b_1 \\ 0 & 0 & 0 & \cdots & 0 & b_1 & b_2  \\ 0 & 0 & 0 & \cdots & b_1 & b_2 & b_3 \\ \vdots & \vdots & \vdots &  & \vdots & \vdots & \vdots \\ b_1 & b_2 & b_3 & \cdots & b_{p-2} & b_{p-1} & b_p  \end{array}\right].$$
There are unique real numbers $k_1, \cdots, k_p$ such that if we set
$$K = \left[\begin{array}{cccc} k_p & k_{p-1} & \cdots & k_1 \\ 0 & k_p & \cdots &k_2 \\ 
\vdots & \vdots & & \vdots \\ 0 & 0 & \cdots & k_p \end{array}\right],$$
then $K^{\textsf{T}} T_{\alpha}^{\textsf{T}}JT_{\alpha} K$ becomes
$$K^{\textsf{T}} T_{\alpha}^{\textsf{T}}JT_{\alpha} K = \left[\begin{array}{ccccc} 0 & 0 & \cdots & 0 & b_1 \\ 0 & 0 & \cdots & b_1 & 0 \\ 
\vdots & \vdots & & \vdots & \vdots \\ 0 & b_1 & \cdots & 0 & 0
\\ b_1 & 0 & \cdots & 0 & 0\end{array}\right].$$
If $\alpha'$ is a cycle in $\mathcal{K}(A)$ whose terminal vector is $w =  \sum_{i=1}^p \, k_i u_i$,  %$$\alpha'=\{A^{p-1}w, A^{p-2}w, \cdots w \},$$ 
then we have $|\alpha'|=p$ and  
$$\langle A^{i}w, A^{j}w \rangle_J = \begin{cases} b_1 \qquad \text{if } j=p-1-i \\ 0 \; \qquad \text{otherwise} \end{cases}$$ 
for each $ 0 \leq i, j \leq p-1$.
%We can replace $\alpha$ with $\alpha_0$. 
If we replace $\alpha$ with $\alpha'$ for each $\alpha$ in $\mathcal{E}(A)$, then the result follows. 

(3) Suppose that $\mathcal{E}(A)$ has properties (1) and (2) and that $\mathcal{E}(A)$ is the union of disjoint cycles $\alpha_1, \cdots, \alpha_r$ of generalized eigenvectors of $A$ corresponding to $0$ for some $r>1$ with $|\alpha_1| \leq |\alpha_2| \leq \cdots \leq |\alpha_r|$. Assuming that
$$\text{span}(\alpha_i) \perp_J \text{span}(\alpha_j) \qquad (i, j =1, \cdots, r-1; i \neq j),$$
we will construct a cycle $\beta$ such that 
the union of the cycles $\alpha_1, \cdots, \alpha_{r-1}, \beta$ forms a basis for $\mathcal{K}(A)$ and
that $\alpha_i$ is orthogonal to $\beta$ with respect to $J$ for each $i=1, \cdots, r-1$.

Suppose that $\alpha = \{u_1, \cdots, u_p\}$ is a cycle in $\mathcal{E}(A)$ which is distinct from 
$\alpha_r = \{w_1, \cdots, w_q\}$. 
We set
$$\langle u_1, u_p\rangle_J = a \, (\neq 0), \qquad \langle u_i, w_{q} \rangle_J = b_i \qquad (i=1, \cdots, p)$$
and
$$z = w_q - \frac{b_1}{a} u_{p} - \frac{b_2}{a} u_{p-1} - \cdots - \frac{b_p}{a}{u_1}.$$
Let $\beta$ denote the cycle whose terminal vector is $z$.  

We first show that 
$u_1 \perp_J \text{span}(\beta)$.
Direct computation yields 
\begin{equation}\label{inductively}
\langle u_1, z \rangle_J = 0.
\end{equation}
Since $Au_1 = 0$, it follows that
$$\langle u_1, A^jz \rangle_J = 0 \qquad (j=1, \cdots, q-1)$$
by (\ref{6}).
Thus, $\langle u_1, A^jz \rangle_J = 0$ for all $j=1, \cdots, q$.

Now, we show that 
$u_2 \perp_J \text{span}(\beta)$.
Direct computation yields
$$\langle u_2, z \rangle_J = 0.$$
From $A^2u_2 = 0$, it follows that
$$\langle u_2, A^jz \rangle_J = 0 \qquad (j=2, \cdots, q-1).$$
It remains to show that $\langle u_2, Az \rangle_J=0$ but this is an immediate consequence of (\ref{6}) and (\ref{inductively}). 

Applying this process to each $u_i$ inductively, the result follows. \hfill $\Box$
\end{pf}

\begin{cor}
There is a basis $\mathcal{E}(A) \in \mathcal{B}as(\mathcal{K}(A))$ such that 
if $u$ is the terminal vector of a cycle $\alpha$ in $\mathcal{E}(A)$ with $|\alpha|=p$, then  $v=A^{p-1-k}u$ is the unique vector in $\mathcal{E}(A)$ satisfying
$$\langle A^ku, v \rangle_J \neq 0$$
for each $k = 0, 1, \cdots p-1$.
\end{cor}

In the rest of the section, 
we investigate a relationship between $\langle \;\; , \;\; \rangle_J$ and $\langle \;\; , \;\; \rangle_K$ on bases $\mathcal{E}(A) \in \mathcal{B}as(\mathcal{K}(A))$ and $\mathcal{E}(B) \in \mathcal{B}as(\mathcal{K}(B))$
when there is a $D_{\infty}$-HEE between two flip pairs $(A, J)$ and $(B, K)$. 
Throughout the section, we assume $(A, J)$ and $(B, K)$ are flip pairs with $|\mathcal{B}_1(\textsf{X}_A)|=m$ and $|\mathcal{B}_1(\textsf{X}_B)|=n$ and $(D, E)$ is a $D_{\infty}$-HEE from $(A, J)$ to $(B, K)$.

We note that $E=KD^{\textsf{T}}J$ implies 
$$\langle u, Dv \rangle_J = \langle Eu, v \rangle_K \qquad (u \in \mathbb{R}^m, v \in \mathbb{R}^n). $$
From this, we see that 
$\text{Ker}(E)$ and $\text{Ran}(D)$ are mutually orthogonal with respect to $J$ and that
$\text{Ker}(D)$ and $\text{Ran}(E)$ are mutually orthogonal with respect to $K$, that is,
\begin{equation} \label{eq10}
\text{Ker}(E) \perp_J \text{Ran}(D) \qquad \text{and} \qquad \text{Ker}(D) \perp_K \text{Ran}(E).
\end{equation}

\begin{lemma}\label{connections}
There exist bases $\mathcal{E}(A) \in \mathcal{B}as(\mathcal{K}(A))$ and $\mathcal{E}(B) \in \mathcal{B}as(\mathcal{K}(B))$ having the following properties.
\newline
\emph{(1)} Suppose that $\alpha$ is a cycle in $\mathcal{E}(A)$ with $|\alpha|=p$ and $u =\emph{Ter}(\alpha)$. Then we have
\begin{equation}\label{eq 3.6}
u \in \emph{Ran}(D) \qquad \Leftrightarrow \qquad A^{p-1}u \notin \emph{Ker}(E)
\end{equation}
\emph{(2)} Suppose that $\beta$ is a cycle in $\mathcal{E}(B)$ with $|\beta|=p$ and $v =\emph{Ter}(\beta)$. Then we have
$$v \in \emph{Ran}(E) \qquad \Leftrightarrow \qquad B^{p-1}v \notin \emph{Ker}(D).$$
\end{lemma}

\begin{pf}
We only prove (\ref{eq 3.6}).
Suppose that $\mathcal{E}(A) \in \mathcal{B}as(\mathcal{K}(A))$ has properties (1), (2) and (3) from Lemma \ref{generalized eigenvectors}.
Since $\langle A^{p-1}u, u \rangle_J \neq 0$, (\ref{eq 3.6}) follows from (\ref{eq10}).
 
Suppose that $u \notin \text{Ran}(D)$.
To draw contradiction, we assume that $A^{p-1}u \notin \text{Ker}(E)$. 
By non-degeneracy of $\langle \;\; , \;\; \rangle_K$, there is a $v \in \mathcal{K}(B)$ such that
$\langle EA^{p-1}u, v \rangle_K \neq 0$, or equivalently, $\langle A^{p-1}u, Dv \rangle_J \neq 0$.
This is a contradiction because 
$\langle A^{p-1}u, u \rangle_J \neq 0$ and $\langle A^{p-1}u, w \rangle_J = 0$ for all $w \in \mathcal{E}(A) \setminus \{u\}$.
\hfill $\Box$
\end{pf}

Now we are ready to prove Proposition B. We first indicate some notation.  
When $p \in  \mathcal{I}nd(\mathcal{K}(A))$,
let $\mathcal{E}_p(A; \partial_{D, E}^-)$ denote the union of cycles $\alpha$ in $\mathcal{E}_p(A)$ such that 
$\text{Ter}(\alpha) \notin \text{Ran}(D)$ 
and let $\mathcal{E}_p(A; \partial_{D, E}^+)$ denote the union of cycles $\alpha$ in $\mathcal{E}_p(A)$ such that $\text{Ter}(\alpha) \in \text{Ran}(D)$.
With this notation, Proposition B can be rewritten as follows.

\begin{thmB} 
If $(D, E) : (A, J) \threesim (B, K)$, then there exist bases $\mathcal{E}(A) \in \mathcal{B}as(\mathcal{K}(A))$ and $\mathcal{E}(B) \in \mathcal{B}as(\mathcal{K}(B))$ having the following properties.
\newline
\emph{(1)} Suppose that $p \in  \mathcal{I}nd(\mathcal{K}(A))$ and $\alpha$ is a cycle in $\mathcal{E}_p(A; \partial_{D, E}^+)$ with $\text{Ter}(\alpha)=u$. 
There is a cycle $\beta$ in $\mathcal{E}_{p+1}(B; \partial_{E, D}^-)$ such that if $\text{Ter}(\beta) = v$, then $Dv=u$.
In this case, we have 
\begin{equation}\label{sym2}
\langle A^{p-1}u, u \rangle_J =  \langle B^{p}v, v \rangle_K.
\end{equation} 
\emph{(2)} Suppose that $p \in  \mathcal{I}nd(\mathcal{K}(A))$, $p>1$ and $\alpha$ is a cycle in $\mathcal{E}_p(A; \partial_{D, E}^-)$ with $\text{Ter}(\alpha) = u$. 
There is a cycle $\beta$ in $\mathcal{E}_{p-1}(B; \partial_{E, D}^+)$ such that if $\text{Ter}(\beta) = v$, then $v=Eu$.
In this case, we have 
\begin{equation}\label{sym1}
\langle A^{p-1}u, u \rangle_J =\langle B^{p-2}v, v \rangle_K.
\end{equation}
\end{thmB}

\begin{pf}
If we define zero-one matrices $M$ and $F$ by
$$M =\left[\begin{array}{cc} 0 & D \\ E & 0  \end{array}\right] \qquad \text{and} \qquad F = \left[\begin{array}{cc} J & 0 \\ 0 & K  \end{array}\right],$$
then $(M, F)$ is a flip pair. 
Suppose that $\mathcal{E}(A) \in \mathcal{B}as(\mathcal{K}(A))$ and $\mathcal{E}(B) \in \mathcal{B}as(\mathcal{K}(B))$ have  properties (1), (2) and (3) from Lemma \ref{generalized eigenvectors}.
If we set
$$\mathcal{E}(A) \oplus 0^n = \bigg\{\left[\begin{array}{c} u \\ 0  \end{array}\right] : u \in \mathcal{E}(A) \text{ and } 0 \in \mathbb{R}^n\bigg\}$$
and
$$0^m \oplus \mathcal{E}(B) = \bigg\{\left[\begin{array}{c} 0 \\ v  \end{array}\right] : v \in \mathcal{E}(B) \text{ and } 0 \in \mathbb{R}^m\bigg\},$$ 
then the elements in $\mathcal{E}(A) \oplus 0^n$ or $0^m \oplus \mathcal{E}(B)$ belong to $\mathcal{K}(M)$. 
Conversely, every vector in $\mathcal{K}(M)$ can be expressed as linear combination of vectors in $\mathcal{E}(A) \oplus 0^n$ and $0^m \oplus \mathcal{E}(B)$.
Thus, the set
$\mathcal{E}(M) = \{\mathcal{E}(A) \oplus 0^n\} \cup \{0^m \oplus \mathcal{E}(B)\}$
becomes a basis for $\mathcal{K}(M)$.

If $\alpha$ is a cycle in $\mathcal{E}(M)$, then $|\alpha|$ is an odd number by Lemma \ref{connections}.
If $|\alpha|=2p-1$ for some positive integer $p$, 
then $\alpha$ is one of the following forms:
$$\bigg\{\left[\begin{array}{c} A^{p-1}u \\ 0  \end{array}\right], \left[\begin{array}{c} 0 \\ B^{p-2}Eu  \end{array}\right], \left[\begin{array}{c} A^{p-2}u \\ 0  \end{array}\right],  \cdots,  \left[\begin{array}{c} Au \\ 0  \end{array}\right], \left[\begin{array}{c} 0 \\ Eu  \end{array}\right], \left[\begin{array}{c} u \\ 0  \end{array}\right]\bigg\}$$
or
$$\bigg\{\left[\begin{array}{c} 0 \\ B^{p-1}v  \end{array}\right], \left[\begin{array}{c} A^{p-2}Dv \\ 0   \end{array}\right], \left[\begin{array}{c}  0 \\ B^{p-2}v   \end{array}\right],  \cdots,  \left[\begin{array}{c} 0 \\ Bv  \end{array}\right], \left[\begin{array}{c} Dv \\ 0  \end{array}\right], \left[\begin{array}{c} 0 \\ v  \end{array}\right]\bigg\}.$$
The formulas (\ref{sym2}) and (\ref{sym1}) are followed from (\ref{zerovalue}).
\hfill $\Box$
\end{pf}

Suppose that $\mathcal{E}(A) \in \mathcal{B}as(\mathcal{K}(A))$ has property (1) from Lemma \ref{generalized eigenvectors}.
If $\alpha$ is a cycle in $\mathcal{E}(A)$, we define the sign of $\alpha$ by
$$\text{sgn}(\alpha) = \begin{cases} +1 \qquad \text{if } \langle \text{Ini}(\alpha), \text{Ter}(\alpha) \rangle_J >0 \\ -1 \qquad \text{if } \langle \text{Ini}(\alpha), \text{Ter}(\alpha) \rangle_J <0. \end{cases}$$
We define the sign of $\mathcal{E}_p(A)$ for each $p \in  \mathcal{I}nd(\mathcal{K}(A))$ by
$$\text{sgn}(\mathcal{E}_p(A)) = \prod_{\{\alpha : \alpha \text{ is a cycle in } \mathcal{E}_p(A)\}} \text{sgn}(\alpha).$$
When $(D, E): (A, J) \threesim (B, K)$, we define the signs of $\mathcal{E}_p(A; \partial_{D, E}^+)$ and $\mathcal{E}_p(A; \partial_{D, E}^-)$ for each $p \in  \mathcal{I}nd(\mathcal{K}(A))$ in similar ways.

Proposition B says that if $(D, E): (A, J) \threesim (B, K)$, there exist bases $\mathcal{E}(A) \in \mathcal{B}as(\mathcal{K}(A))$ and $\mathcal{E}(B) \in \mathcal{B}as(\mathcal{K}(B))$ such that
$$\text{sgn}(\mathcal{E}_p(A; \partial_{D, E}^+)) = \text{sgn}(\mathcal{E}_{p+1}(B; \partial_{E, D}^-))) \qquad (p \in  \mathcal{I}nd(\mathcal{K}(A))),$$
and 
$$\text{sgn}(\mathcal{E}_p(A; \partial_{D, E}^-)) = \text{sgn}(\mathcal{E}_{p-1}(B; \partial_{E, D}^+)) \qquad (p \in  \mathcal{I}nd(\mathcal{K}(A)) ; p>1).$$
In Proposition \ref{prop: cycle of length one} below, we will see that the sign of $\mathcal{E}_1(A; \partial_{D, E}^-)$ is always $+1$ if $\mathcal{E}_1(A; \partial_{D, E}^-)$ is non-empty. 
We first prove Proposition C. 

\begin{pfC}
Suppose that $\mathcal{E}(A) \in \mathcal{B}as(\mathcal{K}(A))$ has properties (1), (2) and (3) from Lemma \ref{generalized eigenvectors} and that $p \in  \mathcal{I}nd(\mathcal{K}(A))$.
We denote the terminal vectors of the cycles in $\mathcal{E}_p(A)$ by $u_{(1)}, \cdots, u_{(q)}$.
Suppose that $P$ is the $m \times q$ matrix whose $i$-th column is $u_{(i)}$ for each $i=1, \cdots, q$.
If we set $M=(A^{p-1}P)^{\textsf{T}} J  P$, then the entry of $M$ is given by
$$M(i, j) = \begin{cases} \langle A^{p-1}u_{(i)}, u_{(i)} \rangle_J \quad \text{if } i=j \\ 0 \qquad \qquad \qquad \quad \;\; \text{otherwise} \end{cases}$$
and the sign of $\mathcal{E}_p(A)$ is determined by the product of the diagonal entries of $M$, that is, 
$$\text{sgn}(\mathcal{E}_p(A)) = \begin{cases} +1 \qquad \text{if } \prod_{i=1}^q M(i, i)>0  \\ -1 \qquad \text{if } \prod_{i=1}^q M(i, i)<0.\end{cases}$$

Suppose that $\mathcal{E}'(A) \in \mathcal{B}as(\mathcal{K}(A))$ is another basis having property (1) from Lemma \ref{generalized eigenvectors}.
Then obviously $\mathcal{E}'_p(A)$ is the union of $q$ disjoint cycles.
If $w$ is the terminal vector of a cycle in $\mathcal{E}'_p(A)$, then $w$ can be expressed as a linear combination of vectors in $\mathcal{E}(A) \cap \text{Ker}(A^p)$, that is,
$$w = \sum_{\substack{k \in \mathbb{R}\\ u \in \mathcal{E}(A) \cap \text{Ker}(A^p)}} ku.$$
If $u \in \mathcal{E}_k(A)$ for $k < p$, then $A^{p-1} u = 0$.
If $u \in \mathcal{E}_k(A)$ for $k>p$ or $u \in \mathcal{E}_p(A)$ and $u$ is not a terminal vector, then $\langle A^{p-1}u, u \rangle_J = 0$ by property (2) from Lemma \ref{generalized eigenvectors}.
This means that the sign of $\mathcal{E}'_p(A)$ is not affected by vectors $u \in \mathcal{E}_k(A)$ for $k \neq p$ or $u \in \mathcal{E}_p(A) \setminus \text{Ter}(\mathcal{E}_p(A))$.
In other words, if we write
$$w = \sum_{i=1}^q \, k_i u_{(i)} + \sum_{u \in u \in \mathcal{E}(A) \setminus \text{Ter}(\mathcal{E}_p(A))} ku \qquad (k_i, k \in \mathbb{R}),$$
then we have
$$\langle A^{p-1}w, w \rangle_J = \langle A^{p-1} \sum_{i=1}^q \, k_i u_{(i)} , \sum_{i=1}^q \, k_i u_{(i)} \rangle_J.$$
To compute the sign of $\mathcal{E}'_p(A)$, we may assume that
$$w =  \sum_{i=1}^q \, k_i u_{(i)} \qquad(k_1, \cdots, k_q \in \mathbb{R}).$$   
We denote the terminal vectors of the cycles in $\mathcal{E}'(A)$ by $w_{(1)}, \cdots, w_{(q)}$ and let 
$Q$ be the $m \times q$ matrix whose $i$-th column is $w_{(i)}$ for each $i=1, \cdots, q$.
If we set $N=(A^{p-1}Q)^{\textsf{T}} J  Q$, then $\prod_{i=1}^q N(i, i) \neq 0$ since $\mathcal{E}'(A)$ has property (1) from Lemma \ref{generalized eigenvectors}.
So we have 
$$\text{sgn}(\mathcal{E}_p'(A)) = \begin{cases} +1 \qquad \text{if } \prod_{i=1}^q N(i, i)>0  \\ -1 \qquad \text{if } \prod_{i=1}^q N(i, i)<0.\end{cases}$$
It is obvious that there is a non-singular matrix $R$ such that $PR=Q$. Since $N=R^{\textsf{T}}MR$ and $M$ is a diagonal matrix, it follows that
$$\prod_{i=1}^q \, M(i, i) > 0 \qquad \Leftrightarrow \qquad \prod_{i=1}^q \, N(i, i)> 0$$
and
$$\prod_{i=1}^q \, M(i, i) < 0 \qquad \Leftrightarrow \qquad \prod_{i=1}^q \, N(i, i) < 0.$$
\hfill $\Box$
\end{pfC}

\begin{prop}\label{prop: cycle of length one}
Suppose that $(D, E): (A, J) \threesim (B, K)$ and that $\mathcal{I}nd(\mathcal{K}(A))$ contains $1$. There is a basis $\mathcal{E}(A) \in \mathcal{B}as(\mathcal{K}(A))$ such that
if $\alpha$ is a cycle in $\mathcal{E}_1(A; \partial_{D, E}^-)$, then $\emph{sgn}(\alpha)=+1$. 
Hence, we have 
$$\emph{sgn}(\mathcal{E}_1(A; \partial_{D, E}^-))=+1.$$
\end{prop}

\begin{pf}
Suppose that $\mathcal{U}$ is a basis for the subspace $\text{Ker}(A)$ of $\mathcal{K}(A)$. 
We may assume that for each $u \in \mathcal{U}$,
\begin{equation}\label{assumption}
a_1, a_2 \in \mathcal{B}_1(\textsf{X}_A), u(a_1) \neq 0 \text{ and } \mathcal{P}_E(a_1) \cap \mathcal{P}_E(a_2) = \varnothing 
\quad \Rightarrow \quad u(a_2)=0
\end{equation}
for the following reason.
If $u(a_2) \neq 0$, then we define
$u_1$ and $u_2$ by
$$u_1(a) = \begin{cases} u(a) \qquad \text{if } \mathcal{P}_E(a_1) \cap \mathcal{P}_E(a) \neq \varnothing \\ 0 \qquad \quad \; \text{otherwise} \end{cases}$$
and
$$u_2(a) = \begin{cases} u(a) \qquad \text{if } \mathcal{P}_E(a_2) \cap \mathcal{P}_E(a) \neq \varnothing \\ 0 \qquad \quad \; \text{otherwise} \end{cases}.$$
It is obvious that $\{u_1, u_2\}$ is linearly independent.
If we set $u_3 = u - u_1 -u_2$ and $u_3 \neq 0$, then obviously, $\{u_1, u_2, u_3\}$ is also linearly independent.
We set
$$\mathcal{U}' =\mathcal{U} \cup  \{u_1, u_2, u_3\}  \setminus\{u\}.$$ 
If necessary, we apply the same process to $u_3$ and to each $u \in \mathcal{U}$ so that every element in $\mathcal{U}'$ satisfies (\ref{assumption}) and then we remove some elements in $\mathcal{U}'$ so that it becomes a basis for $\text{Ker}(A)$.

We first show the following:
$$u \in \mathcal{U} \qquad \Rightarrow \qquad u(\tau_J(a))u(a) \geq 0 \qquad \forall \, a \in \mathcal{B}_1(\textsf{X}_A).
$$
Suppose that $u \in \mathcal{U}$, $a_0 \in \mathcal{B}_1(\textsf{X}_A)$ and that $u(a_0) \neq 0$.
If $a_0 = \tau_J(a_0)$, then $u(\tau_J(a_0))u(a_0) > 0$ and we are done.
When $a_0 \neq \tau_J(a_0)$ and $u(\tau_J(a_0)) = 0$, there is nothing to do.
So we assume $a_0 \neq \tau_J(a_0)$ and $u(\tau_J(a_0)) \neq 0$.
If there were $b \in \mathcal{P}_E(a_0) \cap \mathcal{P}_E(\tau_J(a_0))$,
then we would have
$$1 \geq B(b, \tau_K(b)) \geq E(b, a_0) D(a_0, \tau_K(b)) + E(b, \tau_J(a_0))D(\tau_J(a_0), \tau_K(b))=2$$
from $E=KD^{\textsf{T}}J$.
Thus, we have $\mathcal{P}_E(a_0) \cap \mathcal{P}_E(\tau_J(a_0)) = \varnothing$ and
this implies $u(\tau_J(a_0))=0$ by assumption (\ref{assumption}).

Now, we denote the intersection of $\mathcal{U}$ and $\mathcal{E}_1(A; \partial_{D, E}^-)$ by $\mathcal{V}$ and assume that the elements of $\mathcal{V}$ are $u_1$, $\cdots$, $u_k$, that is, 
$$\mathcal{V} = \mathcal{U} \cap \mathcal{E}_1(A; \partial_{D, E}^-) = \{u_1, \cdots, u_k \}.$$
By Lemma \ref{non-singular} and (\ref{eq10}), for each $u \in \mathcal{V}$ there is a $v \in \mathcal{V}$ such that $\langle u, v \rangle_J \neq 0$. 
If  $\langle u_1, u_1 \rangle_J = 0$, we choose $u_i \in \mathcal{V}$ such that
$\langle u_1, u_i \rangle_J \neq 0$.
There are real numbers $k_1, k_2$ such that
$\{u_1+k_1u_i, u_1+k_2u_i\}$ is linearly independent and that both $\langle u_1+k_1u_i, u_1+k_1u_i\rangle_J$ and $\langle u_1+k_2u_i, u_1+k_2u_i \rangle_J$ are positive.
We replace $u_1$ and $u_i$ with $u_1+k_1u_i$ and $u_1+k_2u_i$.
Continuing this process, we can construct a new basis for $\mathcal{E}_1(A; \partial_{D, E}^-)$ such that
if $\alpha$ is a cycle in $\mathcal{E}_1(A; \partial_{D, E}^-)$, then
$\text{sgn}(\alpha)=+1$. 
\hfill $\Box$
\end{pf}

Suppose that $\mathcal{E}(A) \in \mathcal{B}as(\mathcal{K}(A))$ has property (1) from Lemma \ref{generalized eigenvectors}. 
We arrange the elements of $ \mathcal{I}nd(\mathcal{K}(A)) = \{p_1, p_2, \cdots, p_A\}$  
to satisfy 
$$p_1 < p_2 < \cdots p_A.$$
and write $$\varepsilon_{p}=\text{sgn}(\mathcal{E}_p(A)).$$
If $| \mathcal{I}nd(\mathcal{K}(A))|$=k, then the $k$-tuple $(\varepsilon_{p_1}, \varepsilon_{p_2}, \cdots, \varepsilon_{p_A})$ is called the \textit{flip signature of $(A, J)$} and $\varepsilon_{p_A}$ is called the \textit{leading signature of $(A, J)$}.
The flip signature of $(A, J)$ is denoted by 
$$\emph{F.Sig}(A, J) = (\varepsilon_{p_1}, \varepsilon_{p_2}, \cdots, \varepsilon_{p_A}).$$
When the eventual kernel $\mathcal{K}(A)$ of $A$ is trivial, we write
$$\mathcal{I}nd(\mathcal{K}(A)) = \{0\}$$ 
and define the flip signature of $(A, J)$ by
$$\emph{F.Sig}(A, J) = (+1).$$
We have seen that both the flip signature and the leading signature are independent of the choice of basis $\mathcal{E}_A \in \mathcal{B}as(\mathcal{K}(A))$ as long as $\mathcal{E}_A$ has property (1) from Lemma \ref{generalized eigenvectors}.

In the next section, we prove Proposition A and in Section 5, we prove Theorem D.

\section{Proof of Proposition A} \label{proof of thm1}
We start with the notion of $D_{\infty}$-higher block codes. (See \cite{Ki, LM} for more details about higher block codes.)
We need some notation.
Suppose that $(X, \sigma_X)$ is a shift space over a finite set $\mathcal{A}$ and that $\varphi_{\tau}$ is a one-block flip for $(X, \sigma_X)$ defined by
$$\varphi_{\tau}(x)_i = \tau(x_{-i}) \qquad (x \in X; i \in \mathbb{Z}).$$ 
For each positive integer $n$, we define the $n$-initial map $i_n : \bigcup_{k=n}^{\infty} \mathcal{B}_k(X) \rightarrow \mathcal{B}_n(X)$, the $n$-terminal map $t_n : \bigcup_{k=n}^{\infty} \mathcal{B}_k(X) \rightarrow \mathcal{B}_n(X)$ and the mirror map $\mathcal{M}_n : \mathcal{A}^n \rightarrow \mathcal{A}^n$ by
$$i_n(a_1 a_2 \cdots a_m) = a_1 a_2 \cdots a_n \qquad (a_1 \cdots a_m \in \mathcal{B}_m(X); \; m \geq n),$$
$$t_n(a_1 a_2 \cdots a_m) = a_{m-n+1} a_{m-n+2} \cdots a_m \qquad (a_1 \cdots a_m \in \mathcal{B}_m(X); \; m \geq n)$$
and
$$\mathcal{M}_n(a_1 a_2 \cdots a_n) = a_n \cdots a_1 \qquad (a_1 \cdots a_n \in \mathcal{A}^n).$$
For each positive integer $n$, we denote the map
$$a_1 a_2 \cdots a_n \mapsto \tau(a_1) \tau(a_2) \cdots \tau(a_n) \qquad (a_1 \cdots a_n \in \mathcal{A}^n)$$
by $\tau_n : \mathcal{A}^n \rightarrow \mathcal{A}^n$. 
It is obvious that
the restriction of the map $\mathcal{M}_n \circ \tau_n$ to $\mathcal{B}_n(X)$ is a permutation of order $2$.

For each positive integer $n$, we define the $n$-th higher block code
$h_n: X \rightarrow \mathcal{B}_n(X)^{\mathbb{Z}}$ by
$$h_n(x)_ i = x_{[i, i+n-1]} \qquad (x\in X ; i \in \mathbb{Z}).$$
We denote the image of $(X, \sigma_X)$ under $h_n$ by $(X_n, \sigma_n)$ and call $(X_n, \sigma_n)$ the $n$-th higher block shift of $(X, \sigma_X)$. 
If we write $\upsilon = \mathcal{M}_n \circ \tau_n$, then the map $\varphi_{\upsilon}: X_n \rightarrow X_n$ defined by
$$\varphi_{\upsilon}(x)_i = \upsilon(x_{-i}) \qquad (x \in X_n; i \in \mathbb{Z})$$
becomes a natural one-block flip for $(X_n, \sigma_n)$. 
It is obvious that the $n$-th higher block code $h_n$ is a $D_{\infty}$-conjugacy from $(X, \sigma_X, \varphi_{\tau})$ to $(X_n, \sigma_n, (\sigma_n)^{n-1} \circ \varphi_{\upsilon})$. 
We call the $D_{\infty}$-system $(X_n, \sigma_n, \varphi_{\upsilon})$ the \textit{$n$-th higher block $D_{\infty}$-system} of $(X, \sigma_X, \varphi_{\tau})$.

For notational simplicity, we drop the subscript $n$ and write $\tau=\tau_n$ if the domain of $\tau_n$ is clear in the context.

Suppose that $(A, J)$ is a flip pair. Then the flip pair $(A_n, J_n)$ for the $n$-th higher block $D_{\infty}$-system $(\textsf{X}_n, \sigma_n, \varphi_n)$ of $(\textsf{X}_A, \sigma_A, \varphi_{A, J})$ consists of $\mathcal{B}_n(\textsf{X}_A) \times \mathcal{B}_n(\textsf{X}_A)$ zero-one matrices $A_n$ and $J_n$ defined by  
$$A_n(u, v) = \begin{cases} 1 \qquad \text{if } t_{n-1}(u)=i_{n-1}(v),\\0 \qquad \text{otherwise}\end{cases} \qquad \big(u, v \in \mathcal{B}_n(\textsf{X}_A)\big)$$
and
$$J_n(u, v) = \begin{cases} 1 \qquad \text{if } v = (\mathcal{M} \circ 
\tau_J)(u),\\0 \qquad \text{otherwise}\end{cases} \qquad \big(u, v \in \mathcal{B}_n(\textsf{X}_A)\big).$$

In the following lemma, we prove that there is a $D_{\infty}$-SSE from $(A, J)$ to $(A_n, J_n)$.
\begin{lemma}
\label{lemma: 2.3}
If $n$ is a positive integer greater than $1$, then we have
$$(A_1, J_1) \approx (A_{n}, J_{n})\; (\emph{lag} \; n-1).$$
\end{lemma}

\begin{pf}
For each $k=1, 2, \cdots, n-1$,
we define a zero-one $\mathcal{B}_k(\textsf{X}_A) \times \mathcal{B}_{k+1}(\textsf{X}_A)$ matrix $D_k$ and a zero-one $ \mathcal{B}_{k+1}(\textsf{X}_A) \times \mathcal{B}_k (\textsf{X}_A)$ matrix $E_k$ by 
$$D_k(u, v)=\begin{cases} 1 \qquad \text{if} \; u=i_{k}(v), \\ 0 \qquad \text{otherwise}, \end{cases} \qquad \big(u \in \mathcal{B}_k(\textsf{X}_A), \, v \in \mathcal{B}_{k+1}(\textsf{X}_A)\big)$$
and
$$E_k(v, u)=\begin{cases} 1 \qquad \text{if} \; u=t_k(v), \\ 0 \qquad \text{otherwise} \end{cases} \qquad \big(u \in \mathcal{B}_k(\textsf{X}_A), \, v \in \mathcal{B}_{k+1}(\textsf{X}_A)\big).$$
It is straightforward to see that $(D_k, E_k): (A_k, J_k) \threesim (A_{k+1}, J_{k+1})$ for each $k$.
\hfill $\Box$
\end{pf}

In the proof of Lemma \ref{lemma: 2.3},
$(\textsf{X}_{A_{k+1}}$, $\sigma_{A_{k+1}}$, $\varphi_{A_{k+1}, J_{k+1}})$ is equal to the second higher block $D_{\infty}$-system of $(\textsf{X}_{A_k}$, $\sigma_{A_k}$, $\varphi_{A_k, J_k})$ by recoding of symbols and 
the half elementary conjugacy 
$$\gamma_{D_k, E_k} : (\textsf{X}_{A_k}, \sigma_{A_k}, \varphi_{A_k, J_k}) \rightarrow (\textsf{X}_{A_{k+1}}, \sigma_{A_{k+1}}, \sigma_{A_{k+1}} \circ \varphi_{A_{k+1}, J_{k+1}})$$ induced by $(D_k, E_k)$ can be regarded as the second $D_{\infty}$-higher block code
for each $k=1, 2, \cdots, n-1$. 
A $D_{\infty}$-HEE $(D, E):(A, J)\threesim (B, K)$ is said to be \textit{a complete $D_{\infty}$-half elementary equivalence from $(A, J)$ to $(B, K)$} if $\gamma_{D, E}$ is the second $D_{\infty}$-higher block code. 

In the rest of the section, we prove Proposition A.

\begin{pfA}
We only prove (a). One can prove (b) in a similar way.

We denote the flip pairs for the $n$-th higher block $D_{\infty}$-system of $(\textsf{X}_{A}, \sigma_A, \varphi_{A, J})$ by $(A_n, J_n)$ for each positive integer $n$.
If $\psi :(\textsf{X}_{A}, \sigma_A, \varphi_{A, J}) \rightarrow (\textsf{X}_{B}, \sigma_B, \varphi_{B, K})$ is a $D_{\infty}$-conjugacy, then there are nonnegative integers $s$ and $t$ and a block map $\Psi : \mathcal{B}_{s+t+1}(\textsf{X}_A) \rightarrow \mathcal{B}_1(\textsf{X}_B)$ such that
$$\psi(x)_i = \Psi(x_{[i-s, i+t]}) \qquad (x \in \textsf{X}_A; \; i \in \mathbb{Z}).$$
We may assume that $s+t$ is even by extending window size if necessary. 
By Lemma \ref{lemma: 2.3}, there is a $D_{\infty}$-SSE of lag $(s+t)$ from $(A, J)$ to $(A_{s+t+1}, J_{s+t+1})$.
From (\ref{eq: 2.8}), it follows that the ($s+t+1$)-th $D_{\infty}$-higher block code $h_{s+t+1}$
is a $D_{\infty}$-conjugacy. 
It is obvious that there is a $D_{\infty}$-conjugacy $\psi'$ induced by $\psi$ 
satisfying $\psi = \psi' \circ h_{s+t+1}$ and 
$$x, y \in h_{s+t+1}(X) \qquad \text{and} \qquad x_0 = y_0 \qquad \Rightarrow \qquad \psi'(x)_0 = \psi'(y)_0.$$
So we may assume $s=t=0$ and show that there is a $D_{\infty}$-SSE of lag $2l$ from $(A, J)$ to $(B, K)$ for some positive integer $l$. 

If $\psi^{-1}$ is the inverse of $\psi$, there is a nonnegative integer $m$ such that
\begin{equation}\label{eq: 2.10}
y, y' \in \textsf{X}_B \;\; \text{and} \;\; y_{[-m, m]}=y'_{[-m, m]} \quad \Rightarrow \quad \psi^{-1}(y)_0=\psi^{-1}(y')_0
\end{equation}
since $\psi^{-1}$ is uniformly continuous.
For each $k=1$, $\cdots$, $2m+1$, we define a set $\mathcal{A}_k$ by
$$\mathcal{A}_k = \left\{ \left[\begin{array}{c} v \\ w \\ u \end{array} \right]: u, v \in \mathcal{B}_{i}(\textsf{X}_B), w \in \mathcal{B}_{j}(\textsf{X}_A) \text{ and } u \Psi(w) v \in \mathcal{B}_k(\textsf{X}_{B})\right\},$$
where $i=\lfloor \frac{k-1}{2} \rfloor$ and $j=k-2\lfloor\frac{k-1}{2} \rfloor$.
We define $\mathcal{A}_k \times \mathcal{A}_k$ matrices $M_k$ and $F_k$ to be
\begin{eqnarray*}
 M_k \Bigg( \left[\begin{array}{c} v \\ w \\ u \end{array} \right], \left[\begin{array}{c} v' \\ w' \\ u' \end{array} \right] \Bigg) = 1 \quad
&\Leftrightarrow& \quad \left[\begin{array}{c} v \\ \Psi(w) \\ u \end{array} \right] \left[\begin{array}{c} v' \\ \Psi(w') \\ u' \end{array} \right] \in \mathcal{B}_2(\textsf{X}_{B_k}) \\
&&\\
&& \qquad \quad \text{and} \quad ww' \in \mathcal{B}_2(\textsf{X}_{A_j})
\end{eqnarray*}
and
$$
F_k \Bigg( \left[\begin{array}{c} v \\ w \\ u \end{array} \right], \left[\begin{array}{c} v' \\ w' \\ u' \end{array} \right] \Bigg)=1 
\quad \Leftrightarrow \quad \begin{array}{c} \\ u' = (\mathcal{M} \circ \tau_{K}) (v), \; w'= (\mathcal{M} \circ \tau_{J})(w) \\ \\
\qquad \quad \text{and}\;\; v'= (\mathcal{M} \circ \tau_{K})(u) \end{array}
$$
for all $$\left[\begin{array}{c} v \\ w \\ u \end{array} \right] \;, \; \left[\begin{array}{c} v' \\ w' \\ u' \end{array} \right] \in \mathcal{A}_k.$$
A direct computation shows that $(M_k, F_k)$ is a flip pair for each $k$.
Next, we define a zero-one $\mathcal{A}_k \times \mathcal{A}_{k+1}$ matrix $R_k$ and a zero-one $\mathcal{A}_{k+1} \times \mathcal{A}_k$ matrix $S_k$ to be
$$ R_k \Bigg( \left[\begin{array}{c} v \\ w \\ u \end{array} \right], \left[\begin{array}{c} v' \\ w' \\ u' \end{array} \right] \Bigg)  = 1 
\qquad \Leftrightarrow \qquad  \begin{array}{c} \\ u \Psi(w) v = i_k \left( u' \Psi(w') v' \right) \\ \\
\quad \quad \text{and} \;\; t_1(w)=i_1(w') \end{array}$$
and
$$S_k\Bigg( \left[\begin{array}{c} v' \\ w' \\ u' \end{array} \right], \left[\begin{array}{c} v \\ w \\ u \end{array} \right] \Bigg)  =  1 \qquad
\Leftrightarrow \qquad \begin{array}{c} \\ t_k \left( u'\Psi(w')v' \right)=u\Psi(w)v \\ \\
\quad \quad \text{and} \;\; t_1(w')=i_1(w), \end{array}$$
for all 
$$\left[\begin{array}{c} v \\ w \\ u \end{array} \right] \in \mathcal{A}_{k} \qquad \text{and} \qquad \left[\begin{array}{c} v' \\ w' \\ u' \end{array} \right] \in \mathcal{A}_{k+1}.$$
A direct computation shows that
$$(R_k, S_k):(M_k, F_k)\threesim (M_{k+1}, F_{k+1}).$$
Because $M_1=A$ and $F_1=J$, we obtain
\begin{equation}\label{eq: 2.11}
(A, J) \approx (M_{2m+1}, F_{2m+1})\; (\text{lag} \;2m).
\end{equation}

Finally, (\ref{eq: 2.10}) implies that the $D_{\infty}$-TMC determined by the flip pair $(M_{2m+1}$, $F_{2m+1})$ is equal to the $(2m+1)$-th higher block $D_{\infty}$-system of $(\textsf{X}_B$, $\sigma_B$, $\varphi_{K, B})$ by recoding of symbols.
From Lemma \ref{lemma: 2.3}, we have
\begin{equation}\label{eq: 2.12}
(B, K) \approx (M_{2m+1}, F_{2m+1})\; (\text{lag} \;2m).
\end{equation}
From (\ref{eq: 2.11}) and (\ref{eq: 2.12}), it follows that 
$$(A, J) \approx (B, K)\; (\text{lag} \;4m).$$
\hfill $\Box$
\end{pfA}

\section{Proof of Theorem D}\label{sec: fs}

We start with the case where $(B, K)$ in Theorem D is the flip pair for the $n$-th higher block $D_{\infty}$-system of $(\textsf{X}_A, \sigma_A, \varphi_{A, J})$.

\begin{lemma}\label{lem 5.1}
Suppose that $(B, K)$ is the flip pair for the $n$-th higher block $D_{\infty}$-system of $(\textsf{X}_A, \sigma_A, \varphi_{A, J})$.
\newline
\emph{(1)} If $p \in \mathcal{I}nd(\mathcal{K}(A))$, then there is $q \in \mathcal{I}nd(\mathcal{K}((B))$ such that $q = p+n-1$ and that
$$\text{sgn}(\mathcal{E}_p(A)) = \text{sgn}(\mathcal{E}_{q}(B)).$$
\emph{(2)} If $q \in \mathcal{I}nd(\mathcal{K}((B))$ and $q \geq n$, then there is $p \in \mathcal{I}nd(\mathcal{K}(A))$ such that
$q = p + n - 1$ 
and that
$$\text{sgn}(\mathcal{E}_p(A)) = \text{sgn}(\mathcal{E}_q(B)).$$
\emph{(3)} If $q \in \mathcal{I}nd(\mathcal{K}((B))$
and $q < n$, then we have 
$$\text{sgn}(\mathcal{E}_q(B)) = +1.$$ 
\end{lemma}

\begin{pf} 
We only prove the case $n=2$. 
We assume $\mathcal{E}(A) \in \mathcal{B}as(\mathcal{K}(A))$ and $\mathcal{E}(B) \in \mathcal{B}as(\mathcal{K}(B))$ are bases having properties from Proposition B.
Suppose that $\alpha$ is a cycle in $\mathcal{E}_p(A)$ for some $p \in \mathcal{I}nd(\mathcal{K}(A))$ and that $u$ is the initial vector of $\alpha$.
For any $a_1 a_2 \in \mathcal{B}_2(\textsf{X}_A)$, we have
$$Eu\left(\left[\begin{array}{c} a_2 \\ a_1 \end{array}\right]\right) = u(a_2)$$
and this implies that $Eu$ is not identically zero.
By Lemma \ref{connections}, $\alpha$ is a cycle in $\mathcal{E}_p(A; \partial_{D, E}^+)$.
Under the assumption that $\mathcal{E}(A)$ and $\mathcal{E}(B)$ have properties from Proposition B, we can find a cycle $\beta$ in $\mathcal{E}(B)$ such that the initial vector of $\beta$ is $Eu$. 
Thus, we obtain
\begin{equation}\label{eq in Lem 5.1}
\mathcal{E}_p(A; \partial_{D, E}^-) = \varnothing \qquad \text{and} \qquad \mathcal{E}_{p+1}(B; \partial_{E, D}^+) = \varnothing,
\end{equation}
$$p \in  \mathcal{I}nd(\mathcal{K}(A)) \qquad \Leftrightarrow \qquad p+1 \in  \mathcal{I}nd(\mathcal{K}(B)) \qquad (p \geq 1)$$ 
and
$$\text{sgn}(\mathcal{E}_p (A))=\text{sgn}(\mathcal{E}_{p+1}(B)) \qquad (p\in  \mathcal{I}nd(\mathcal{K}(A))).$$
If $\mathcal{E}_1(B) \neq \varnothing$, then $\mathcal{E}_1(B) = \mathcal{E}_1(B; \partial_{E, D}^-)$ by (\ref{eq in Lem 5.1})
and we have
$$\text{sgn}(\mathcal{E}_1(B)) = +1$$
by Proposition \ref{prop: cycle of length one} and Proposition C.
\hfill $\Box$
\end{pf}

\begin{rmk} 
If two $D_{\infty}$-TMCs are finite, then we can directly determine whether or not they are $D_{\infty}$-conjugate.	
In this paper, we do not consider  $D_{\infty}$-TMCs who have finite cardinalities.
Hence, when $(B, K)$ is the flip pair for the $n$-th higher block $D_{\infty}$-system of $(\textsf{X}_A, \sigma_A, \varphi_{A, J})$ for some positive integer $n>1$, $B$ must have zero as its eigenvalue.
\end{rmk}

\begin{pfD}	
Suppose that $(A, J)$ and $(B, K)$ are flip pairs and that $\psi:(\textsf{X}_A, \sigma_A,\varphi_{A, J}) \rightarrow (\textsf{X}_B, \sigma_B, \varphi_{B, K})$ is a $D_{\infty}$-conjugacy.
As we can see in the proof of Proposition A, there is a $D_{\infty}$-SSE from $(A, J)$ to $(B, K)$ consisting of the even number of complete $D_{\infty}$-half elementary equivalences and $(R_k, S_k):(M_k, F_k) \threesim (M_{k+1}, F_{k+1})$ $(k=1, \cdots, 2m)$.
In Lemma \ref{lem 5.1}, we have already seen that Theorem D is true in the case of complete $D_{\infty}$-half elementary equivalences.
So it remains to compare the flip signatures of $(M_k, F_k)$ and $(M_{k+1}, F_{k+1})$ for each $k=1, \cdots, 2m$.
Throughout the proof, we assume $\mathcal{A}_k$ and $(R_k, S_k):(M_k, F_k) \threesim (M_{k+1}, F_{k+1})$ are as in the proof of Proposition A.

We only discuss the following two cases:
\newline
(1) $(R_2, S_2):(M_2, F_2) \threesim (M_3, F_3)$ 
\newline
(2) $(R_3, S_3):(M_3, F_3) \threesim (M_4, F_4)$.
\newline
When $k=1$, $(R_1, S_1)$ is a complete $D_{\infty}$-half elementary conjugacy from $(A, J)$ to $(A_2, J_2)$.
For each $k=4, \cdots, 2m$, one can apply the arguments used in (1) and (2) to $(R_k, S_k):(M_k, F_k) \threesim (M_{k+1}, F_{k+1})$. When $k$ is an even number, the argument used in (1) can be applied and when $k$ is an odd number, the argument used in (2) can be applied.

(1) Suppose that $(B_2, K_2)$ is the flip pair for the second higher block $D_{\infty}$-system of $(\textsf{X}_B, \sigma_B, \varphi_{B, K})$.
We first compare the flip signatures of $(B_2, K_2)$ and $(M_3, F_3)$.
We define a zero-one $\mathcal{B}_2(\textsf{X}_B) \times \mathcal{A}_3$ matrix $U_2$ and a zero-one $\mathcal{A}_3 \times \mathcal{B}_2(\textsf{X}_B)$ matrix $V_2$ by
$$U_2 \Bigg( \left[\begin{array}{c} b_2 \\ b_1 \end{array} \right], \left[\begin{array}{c} d_3 \\ a_2 \\ d_1 \end{array} \right] \Bigg)  = \begin{cases} 1 
\qquad \text{if }  b_1=d_1 \text{ and } \Psi(a_2) = b_2 \\ 0 \qquad \text{otherwise} \end{cases} $$
and
$$V_2\Bigg( \left[\begin{array}{c} d_3 \\ a_2 \\ d_1 \end{array} \right], \left[\begin{array}{c} b_2 \\ b_1 \end{array} \right] \Bigg)  =  \begin{cases} 1 \qquad \text{if } b_2=d_3 \text{ and } \Psi(a_2)=b_1 \\ 0 \qquad \text{otherwise} \end{cases}$$
for all 
$$\left[\begin{array}{c} b_2 \\ b_1 \end{array} \right] \in \mathcal{B}_2(\textsf{X}_B) \qquad \text{and} \qquad \left[\begin{array}{c} d_3 \\ a_2 \\ d_1 \end{array} \right] \in \mathcal{A}_{3}.$$
A direct computation shows that
$$(U_2, V_2):(B_2, K_2)\threesim (M_3, F_3).$$

Remark of Lemma \ref{lem 5.1} says that $\mathcal{K}(B_2)$ is not trivial.
So there is a basis $\mathcal{E}(B_2) \in \mathcal{B}as(\mathcal{K}(B_2))$ for the eventual kernel of $B_2$ having property (1) from Lemma \ref{generalized eigenvectors}. Suppose that 
$\gamma = \{w_1, \cdots, w_p\}$ is a cycle in $\mathcal{E}(B_2)$.
Since
$$V_2w_1 \Bigg( \left[\begin{array}{c} b_3 \\ a_2 \\ b_1 \end{array} \right]\Bigg) = w_1 \Big(\Big[ \begin{array}{c} b_3 \\ b_2 \end{array}\Big]\Big) \qquad \Bigg( \left[\begin{array}{c} b_3 \\ a_2 \\ b_1 \end{array} \right] \in \mathcal{A}_3 \Bigg),$$ 
it follows that $w_1 \notin \text{Ker}(V_2)$.
By Lemma \ref{connections}, $\gamma$ is a cycle in $\mathcal{E}_p(B_2; \partial_{U_2, V_2}^+)$.
Suppose that $\mathcal{E}(M_3) \in \mathcal{B}as(\mathcal{K}(M_3))$ is a basis for the eventual kernel of $M_3$ having property (1) from Lemma \ref{generalized eigenvectors}.
Then it is obvious that for each $p \in  \mathcal{I}nd(\mathcal{K}(B_2))$, we have
\begin{equation}\label{eq in thmD}
\mathcal{E}_p(B_2; \partial_{U_2, V_2}^-) = \varnothing \qquad \text{and} \qquad \mathcal{E}_{p+1}(M_3; \partial_{V_2, U_2}^+) = \varnothing. 
\end{equation}
Hence, 
$$p \in  \mathcal{I}nd(\mathcal{K}(B_2)) \qquad \Leftrightarrow \qquad p+1 \in  \mathcal{I}nd(\mathcal{K}(M_3)) \qquad (p \geq 1)$$ 
and
$$\text{sgn}(\mathcal{E}_p (B_2))=\text{sgn}(\mathcal{E}_{p+1}(M_3)) \qquad (p\in  \mathcal{I}nd(\mathcal{K}(B_2)))$$
by Proposition C.
If $\mathcal{E}_1(M_3) \neq \varnothing$, then $\mathcal{E}_1(M_3) = \mathcal{E}_1(M_3; \partial_{V_2, U_2}^-)$ by (\ref{eq in thmD})
and we have
\begin{equation}\label{0:+1}
\text{sgn}(\mathcal{E}_1(M_3)) = +1
\end{equation}
by Proposition \ref{prop: cycle of length one} and Proposition C.

Now, we compare the flip signatures of $(M_2, F_2)$ and $(M_3, F_3)$.
Let $\beta=\{v_1, \cdots, v_{p+1}\}$ be a cycle in $\mathcal{E}(M_3)$ for some $p \geq 1$.
If $b_1b_2b_3 \in \mathcal{B}_3(\textsf{X}_B)$ and $a_2, a_2' \in \Psi^{-1}(b_2)$, then from $M_3v_2 =v_1$, it follows that
$$v_1 \Bigg( \left[\begin{array}{c} b_3 \\ a_2 \\ b_1 \end{array} \right]\Bigg) = \sum_{a_3 \in \Psi^{-1}(b_3)} \sum_{b_4 \in \mathcal{F}_B(b_3)} v_2 \Bigg( \left[\begin{array}{c} b_4 \\ a_3 \\ b_2 \end{array} \right]\Bigg)$$
and this implies that
$$v_1 \Bigg( \left[\begin{array}{c} b_3 \\ a_2 \\ b_1 \end{array} \right]\Bigg) = v_1 \Bigg( \left[\begin{array}{c} b_3 \\ a_2' \\ b_1 \end{array} \right]\Bigg).$$
Since $v_1$ is a nonzero vector, 
there is a block $b_1b_2b_3 \in \mathcal{B}_3(\textsf{X}_B)$ and a nonzero real number $k$ such that
$$v_1 \Bigg( \left[\begin{array}{c} b_3 \\ a_2 \\ b_1 \end{array} \right]\Bigg) = k \qquad \forall \, a_2 \in \Psi^{-1}(b_2).$$ 
Since $M_3v_1=0$, it follows that
$$ \sum_{a_2 \in \Psi^{-1}(b_2)} \sum_{b_3 \in \mathcal{F}_B(b_2)}  v_1 \Bigg( \left[\begin{array}{c} b_3 \\ a_2 \\ b_1 \end{array} \right]\Bigg) = k \sum_{b_3 \in \mathcal{F}_B(b_2)}  v_1 \Bigg( \left[\begin{array}{c} b_3 \\ a_2 \\ b_1 \end{array} \right]\Bigg)= 0.$$
From this, we see that
$$R_2 v_1 \Big(\Big[ \begin{array}{c} a_2 \\ a_1 \end{array}\Big]\Big) = \sum_{b_3 \in \mathcal{F}_B(b_2)} v_1\Bigg( \left[\begin{array}{c} b_3 \\ a_2 \\ b_1 \end{array} \right]\Bigg)=0$$
for any $a_1 \in \Psi^{-1}(b_1)$ and $a_1 a_2 \in \mathcal{B}_2(\textsf{X}_A)$.
Hence, $v_1 \in \text{Ker}(R_2)$
and $\beta$ is a cycle in $\mathcal{E}_{p+1}(M_3; \partial_{S_2, R_2}^-)$ by Lemma \ref{connections}.
From this, we see that
$$p+1 \in  \mathcal{I}nd(\mathcal{K}(M_3)) \qquad \Leftrightarrow \qquad p \in  \mathcal{I}nd(\mathcal{K}(M_2)) \qquad (p \geq 2)$$
and 
$$2 \in  \mathcal{I}nd(\mathcal{K}(M_3)) \qquad \Leftrightarrow \qquad 1 \in  \mathcal{I}nd(\mathcal{K}(M_2; \partial_{R_2, S_2}^+)).$$
Suppose that $\mathcal{E}(M_2) \in \mathcal{B}as(\mathcal{K}(M_2))$ is a basis for the eventual kernel of $M_2$ having property (1) from Lemma \ref{generalized eigenvectors}. 
If $1 \in \mathcal{I}nd(\mathcal{K}(M_2))$ and $\mathcal{E}_{1}(M_2; \partial_{R_2, S_2}^-)$ is non-empty, then 
we have
$$\text{sgn}(\mathcal{E}_{1}(M_2; \partial_{R_2, S_2}^-)) = +1$$
by Proposition \ref{prop: cycle of length one}, Proposition C and (\ref{eq10}).
Thus, we have
$$\text{sgn}(\mathcal{E}_{p+1} (M_3))=\text{sgn}(\mathcal{E}_{p}(M_2)) \qquad (p+1\in  \mathcal{I}nd(\mathcal{K}(M_3)); p \geq 1).$$
If $1 \in \mathcal{I}nd(\mathcal{K}(M_3))$ and $\mathcal{E}_{1}(M_3; \partial_{S_2, R_2}^+)$ is non-empty, then we have
$$\text{sgn}(\mathcal{E}_{1}(M_3; \partial_{S_2, R_2}^+)) = +1$$
and if $\mathcal{E}_{1}(M_3; \partial_{S_2, R_2}^-)$ is non-empty, then 
we have
$$\text{sgn}(\mathcal{E}_{1}(M_3; \partial_{S_2, R_2}^-)) = +1$$
by (\ref{eq10}), (\ref{0:+1}), Proposition \ref{prop: cycle of length one} and Proposition C.  
As a consequence, the flip signatures of 
$(M_2, F_2)$ and $(M_3, F_3)$ have the same number of $-1$'s and their leading signatures coincide.

(2) Suppose that $\alpha$ is a cycle in $\mathcal{K}(M_3)$ and that $u$ is the initial vector of $\alpha$.
Since
$$S_3u\left(\left[\begin{array}{c} b_4 \\ a_3 \\ a_2 \\ b_1 \end{array}\right]\right) = u \left[\begin{array}{c} b_4 \\ a_3 \\ \Psi(a_2)  \end{array}\right] \qquad \left(\left[\begin{array}{c} b_4 \\ a_3 \\ a_2 \\ b_1 \end{array}\right] \in \mathcal{A}_4\right),$$
it follows that $S_3u$ is not identically zero.
The argument used in the proof of Lemma \ref{lem 5.1} 
completes the proof.
\hfill $\Box$
\end{pfD}

\section{$D_{\infty}$-Shift Equivalence and the Lind Zeta Functions}

We first introduce the notion of $D_{\infty}$-shift equivalence which is an analogue of shift equivalence.
Let $(A, J)$ and $(B, K)$ be flip pairs and let $l$ be a positive integer.
A \textit{$D_{\infty}$-shift equivalence} ($D_{\infty}$-SE) \textit{of lag $l$ from $(A, J)$ to $(B, K)$} is a pair $(D, E)$ of nonnegative integral matrices satisfying
$$A^{l}=DE, \quad B^{l}=ED, \quad AD=DB, \quad \text{and} \quad E=KD^{\textsf{T}}J.$$
We observe that $AD=DB$, $E=KD^{\textsf{T}}J$ and the fact that $(A,J)$ and $(B,K)$ are flip pairs imply $EA=BE$.
If there is a $D_{\infty}$-SE of lag $l$ from $(A, J)$ to $(B, K)$, then we say that $(A, J)$ is $D_{\infty}$-shift equivalent to $(B, K)$ and write %$(D, E): (A, J) \sim (B, K)$ (lag $l$), or simply 
$$(A, J) \sim (B, K) \; (\text{lag } l).$$

Suppose that  
$$(D_1, E_1), (D_2, E_2), \cdots, (D_l, E_l)$$ 
is a $D_{\infty}$-SSE of lag $l$ from $(A, J)$ to $(B, K)$. 
If we set
$$D=D_1 D_2 \cdots D_l \qquad \text{and} \qquad E=E_l \cdots E_2 E_1,$$
then $(D, E)$ is a $D_{\infty}$-SE of lag $l$ from $(A, J)$ to $(B, K)$.
Hence, we have
$$(A, J) \approx (B, K) \;\; (\text{lag } l) \qquad \Rightarrow \qquad (A, J) \sim (B, K) \; (\text{lag } l).$$

In the rest of the section, we review the Lind zeta function of a $D_{\infty}$-TMC.
In \cite{KLP}, an explicit formula for the Lind zeta function of a $D_{\infty}$-system was established. In the case of a  $D_{\infty}$-TMC, the Lind zeta function can be expressed in terms of matrices from flip pairs. We briefly discuss the formula.

Suppose that $G$ is a group and that $\alpha$ is a $G$-action on a set $X$. Let $\mathcal{F}$ denote the set of finite index subgroups of $G$. For each $H \in \mathcal{F}$, we set
$$p_H(\alpha) = |\{ x \in X : \forall \, h \in H \; \alpha(h, x) = x \}|.$$ 
The Lind zeta function $\zeta_{\alpha}$ of the action $\alpha$ is defined by 
\begin{equation}\label{Lind}
\zeta_{\alpha}(t) = \exp \left( \sum_{H \in \mathcal{F}} \frac{p_H(\alpha)}{|G/H|}\, t^{|G/H|}\right).
\end{equation}
It is clear that if $\alpha : \mathbb{Z} \times X \rightarrow X$ is given by $\alpha(n, x) = T^n(x)$, then the Lind zeta function $\zeta_{\alpha}$ becomes the Artin-Mazur zeta function $\zeta_T$ of a topological dynamical system $(X, T)$. 
The formula for the Artin-Mazur zeta function can be found in \cite{AM}. Lind defined the function (\ref{Lind}) in \cite{L} for the case $G = \mathbb{Z}^d$.

Every finite index subgroup of the infinite dihedral group $D_{\infty} =\langle a, b : ab=ba^{-1} \;\, \text{and} \;\, b^2=1\rangle$ can be written in one and only one of the following forms: 
$$\langle a^m \rangle \qquad \text{or} \qquad \langle a^m, a^k b \rangle \qquad (m=1, 2, \cdots\,; k=1, \cdots, m-1)$$
and the index is given by
$$|G_2/ \langle a^m \rangle| = 2m \qquad \text{or} \qquad |G_2/ \langle a^m, a^k b \rangle| = m.$$
Suppose that $(X, T, F)$ is a $D_{\infty}$-system.
If $m$ is a positive integer, then the number of periodic points in $X$ of period $m$ will be denoted by $p_m(T)$:
$$p_m(T) = |\{ x \in X : T^m(x) = x \}|.$$
If $m$ is a positive integer and $n$ is an integer, then $p_{m,n}(T, F)$ will denote the number of points in $X$ fixed by $T^m$ and $T^n \circ F$:
$$p_{m, n}(T, F) = |\{ x \in X : T^m(x) = T^n \circ F (x) = x \}|.$$
Thus, the Lind zeta function $\zeta_{T,F}$ of a $D_{\infty}$-system $(X, T, F)$ is given by
\begin{equation}\label{Lindzeta}
\zeta_{T, F}(t) = \exp \Big(\sum_{m=1}^{\infty} \, \frac{p_m(T)}{2m}t^{2m} +\sum_{m=1}^{\infty} \, \sum_{k=0}^{m-1}\, \frac{p_{m,k}(T, F)}{m}t^{m} \Big).
\end{equation}
It is evident if two $D_{\infty}$-systems $(X, T, F)$ and $(X', T', F')$ are $D_{\infty}$-conjugate, then 
$$p_m(T) = p_m(T') \qquad \text{and} \qquad p_{m, n}(T, F) = p_{m,n}(T', F')$$
for all positive integers $m$ and integers $n$. 
As a consequence, the Lind zeta function is a $D_{\infty}$-conjugacy invariant.

The formula (\ref{Lindzeta}) can be simplified as follows.
Since $T \circ F = F \circ T^{-1}$ and $F^2 = \text{Id}_X,$ it follows that
$$p_{m,n}(T, F) = p_{m,n+m}(T, F) = p_{m,n+2}(T, F)$$
and this implies that
\begin{eqnarray}
&& p_{m,n}(T, F) = p_{m,0}(T, F) \qquad \text{ if } m \text{ is odd}, \label{eq: 4.1}\\ && p_{m,n}(T, F)= p_{m,0}(T, F) \qquad \text{ if } m \text{ and } n \text{ are even}, \nonumber \\ && p_{m,n}(T, F)= p_{m,1}(T, F) \qquad \text{ if } m \text{ is even and } n \text{ is odd}. \nonumber 
\end{eqnarray}
Hence, we obtain
$$\sum_{k=0}^{m-1} \, \frac{p_{m,n}(T, F)}{m} = \begin{cases} \, p_{m,0}(T, F) \qquad \qquad \qquad \qquad \text{if } m \text{ is odd},\\ \\ \, \displaystyle{\frac{p_{m,0}(T, F)+p_{m,1}(T, F)}{2} \qquad \text{if } m \text{ is even}.}\end{cases}$$
Using this,  (\ref{Lindzeta}) becomes
$$\zeta_{\alpha}(t) = {\zeta_T(t^2)}^{1/2} \exp \left( G_{T, F}(t) \right),$$
where
$\zeta_T$ is the Artin-Mazur zeta function of $(X, T)$ and $G_{T, F}$ is given by
$$G_{T, F}(t) = \sum_{m=1}^{\infty} \, \left( p_{2m-1, 0}(T, F) \, t^{2m-1} + \frac{p_{2m, 0}(T, F)+p_{2m, 1}(T, F)}{2} \, t^{2m}\right).$$

If there is a $D_{\infty}$-SSE of lag $2l$ between two flip pairs $(A, J)$ and $(B, K)$ for some positive integer $l$, then $(\textsf{X}_A, \sigma_A, \varphi_{A, J})$ and $(\textsf{X}_B, \sigma_B, \varphi_{B, K})$ have the same Lind zeta function by (1) in Proposition A .
The following proposition says that the Lind zeta function is actually an invariant for $D_{\infty}$-SSE. 
\begin{prop}
\label{prop: 4.1}
If $(X, T, F)$ is a $D_{\infty}$-system, then
$$p_{2m-1, 0}(T, F) = p_{2m-1, 0}(T, T \circ F),$$
$$p_{2m, 0}(T, F) = p_{2m, 1}(T, T \circ F), %\quad \text{and}
$$
$$p_{2m, 1}(T, F) = p_{2m, 0}(T, T \circ F)$$
for all positive integers $m$.
As a consequence, the Lind zeta functions of $(X, T, F)$ and $(X, T, T \circ F)$ are the same.
\end{prop}

\begin{pf}
The last equality is trivially true. To prove the first two equalities, we observe that 
$$T^{m}(x) = F(x)=x \quad \Leftrightarrow \quad T^{m}(Tx)=T\circ(T \circ F)(Tx) = Tx$$
for all positive integers $m$.
Thus, we have
\begin{equation}
\label{eq: 4.2}
p_{m, 0}(T, F)=p_{m, 1}(T, T \circ F) \qquad (m=1, 2, \cdots).
\end{equation}
Replacing $m$ with $2m$ yields the second equality.
From (\ref{eq: 4.1}) and (\ref{eq: 4.2}), the first one follows.
\hfill $\Box$
\end{pf}

When $(A, J)$ is a flip pair, the numbers $p_{m, \delta}(\sigma_A, \varphi_{A, J})$ of fixed points can be expressed in terms of $A$ and $J$ for all positive integers $m$ and $\delta \in \{0, 1\}$.
In order to present it, we indicate notation. 
If $M$ is a square matrix, then $\Delta_M$ will denote the column vector whose $i$-th coordinates are identical with $i$-th diagonal entries of $M$, that is,
$$\Delta_M(i) = M(i, i).$$
For instance, if $I$ is the $2 \times 2$ identity matrix, then 
$$\Delta_I = \left[\begin{array}{c} 1 \\ 1 \end{array}\right].$$
The following proposition is proved in \cite{KLP}.
\begin{prop}
\label{prop: 4.3}
If $(A, J)$ is a flip pair, then
$$p_{2m-1, 0}(\sigma_A, \varphi_{J, A}) = {\Delta_J}^\textsf{T} \left( A^{m-1} \right) \Delta_{AJ},$$
$$p_{2m, 0}(\sigma_A, \varphi_{J, A}) = {\Delta_J}^\textsf{T} \left( A^m \right) \Delta_J,% \quad \text{and}
$$
$$p_{2m, 1}(\sigma_A, \varphi_{J, A}) = {\Delta_{JA}}^\textsf{T} \left( A^{m-1} \right) \Delta_{AJ}$$
for all positive integers $m$.
\end{prop}

\section{Examples}
\label{sec:fifth}
Let $A$ be Ashley's eight-by-eight and let $B$ be the minimal zero-one transition matrix for the full two-shift, that is,
$$A=\left[\begin{array}{rrrrrrrr} 1 & 1 & 0 & 0 & 0 & 0 & 0 & 0 \\ 0 & 0 & 1 & 0 & 0 & 0 & 1 & 0  \\ 0 & 0 & 0 & 1 & 0 & 1 & 0 & 0 \\ 0 & 1 & 0 & 0 & 0 & 0 & 0 & 1 \\ 1 & 0 & 0 & 0 & 1 & 0 & 0 & 0 \\ 0 & 0 & 0 & 0 & 1 & 0 & 0 & 1 \\ 0 & 0 & 1 & 0 & 0 & 1 & 0 & 0 \\ 0 & 0 & 0 & 1 & 0 & 0 & 1 & 0 \end{array}\right] \quad \text{and} \quad B = \left[\begin{array}{rr} 1 & 1 \\ 1 & 1 \end{array}\right].$$
There is a unique one-block flip for $(\textsf{X}_A, \sigma_A)$ and there are exactly two one-block flips for $(\textsf{X}_B, \sigma_B)$.
Those flips are determined by the permutation matrices
$$J=\left[\begin{array}{rrrrrrrr} 0 & 0 & 0 & 0 & 1 & 0 & 0 & 0 \\ 0 & 0 & 0 & 0 & 0 & 1 & 0 & 0  \\ 0 & 0 & 0 & 0 & 0 & 0 & 1 & 0 \\ 0 & 0 & 0 & 0 & 0 & 0 & 0 & 1 \\ 1 & 0 & 0 & 0 & 0 & 0 & 0 & 0 \\ 0 & 1 & 0 & 0 & 0 & 0 & 0 & 0 \\ 0 & 0 & 1 & 0 & 0 & 0 & 0 & 0 \\ 0 & 0 & 0 & 1 & 0 & 0 & 0 & 0 \end{array}\right], 
\quad I = \left[\begin{array}{rr} 1 & 0 \\ 0 & 1 \end{array}\right] \quad \text{and} \quad K = \left[\begin{array}{rr} 0 & 1 \\ 1 & 0 \end{array}\right].$$

\begin{exa}
\label{exa: 1}
Direct computation shows that the number of fixed points of $(\textsf{X}_A, \sigma_A, \varphi_{A, J})$, $(\textsf{X}_B, \sigma_B, \varphi_{B, I})$ and $(\textsf{X}_B, \sigma_B, \varphi_{B, K})$ are as follows:
$$p_m(\sigma_A) = p_m(\sigma_B) = 2^m,$$
$$p_{2m-1, 0}(\sigma_A, \varphi_{A, J}) = p_{2m, 0}(\sigma_A, \varphi_{A, J}) = 0,$$ $$p_{2m, 1}(\sigma_A, \varphi_{A, J}) = \begin{cases} 2^m \quad \text{if } m \neq 6 \\ 80 \quad \; \text{if } m=6, \end{cases}$$
$$p_{2m-1, 0}(\sigma_B, \varphi_{B, I}) = 2^m, \quad p_{2m, 0}(\sigma_B, \varphi_{B, I}) = 2^{m+1}, \quad p_{2m, 1}(\sigma_B, \varphi_{B, I}) = 2^m,$$
$$p_{2m-1, 0}(\sigma_B, \varphi_{B, K}) = p_{2m, 0}(\sigma_B, \varphi_{B, K}) = 0, \quad p_{2m, 1}(\sigma_B, \varphi_{B, K}) = 2^m$$
for all positive integers $m$.
The Lind zeta functions are as follows: 
$$\zeta_{A, J}(t) = \frac{1}{\sqrt{1-2t^2}} \, \exp \left( \frac{t^2}{1-2t^2} +8t^{12}\right),$$
$$\zeta_{B, I}(t) = \frac{1}{\sqrt{1-2t^2}} \, \exp \left( \frac{2t+3t^2}{1-2t^2} \right)$$
and
$$\zeta_{B, K}(t) = \frac{1}{\sqrt{1-2t^2}} \, \exp \left( \frac{t^2}{1-2t^2} \right).$$
As a result, we see that
$$(\textsf{X}_A, \sigma_A, \varphi_{A, J}) \ncong (\text{X}_B, \sigma_B, \varphi_{B, I}),$$
$$(\textsf{X}_A, \sigma_A, \varphi_{A, J}) \ncong (\text{X}_B, \sigma_B, \varphi_{B, K})$$
and
$$(\textsf{X}_A, \sigma_A, \varphi_{B, I}) \ncong (\text{X}_B, \sigma_B, \varphi_{B, K}).$$
\end{exa}
\begin{exa}
In spite of $\zeta_{A, J} \neq \zeta_{B, I}$, $\zeta_{A, J} \neq \zeta_{B, K}$ and $\zeta_{B, I} \neq \zeta_{B, K}$,
there are $D_{\infty}$-SEs between $(A, J)$, $(B, I)$ and $(B, K)$ pairwise.
If $D$ and $E$ are matrices given by
$$D = 2\left[\begin{array}{rr} 1 & 1 \\ 1 & 1 \\ 1 & 1 \\ 1 & 1 \\ 1 & 1 \\ 1 & 1 \\ 1 & 1 \\ 1 & 1 \end{array}\right]\qquad \text{and} \qquad E= 2\left[\begin{array}{rrrrrrrr} 1 & 1 & 1 & 1 & 1 & 1 & 1 & 1 \\ 1 & 1 & 1 & 1 & 1 & 1 & 1 & 1 \end{array}\right],$$
then $(D, E)$ is a $D_{\infty}$-SE of lag $6$ from $(A, J)$ to $(B, K)$ and from $(A, J)$ to $(B, I)$:
$$(D, E): (A, J) \sim (B, I) \; (\text{lag } 6) \quad \text{and} \quad (D, E): (A, J) \sim (B, K) \; (\text{lag } 6).$$ 
%(Actually, $6$ is the smallest lag for a $D_{\infty}$-SE between $(A, J)$ and $(B, K)$.)

Direct computation shows that $(B^l, B^l)$ is
a $D_{\infty}$-SE from $(B, I)$ to $(B, K)$:
$$(B^l, B^l): (B, I) \sim (B, K) \; (\text{lag } 2l)$$
for all positive integers $l$. 
This contrasts with the fact that the existence of SE between two transition matrices implies that the corresponding $\mathbb{Z}$-TMCs share the same Artin-Mazur zeta functions. (See Section 7 in \cite{LM}.)
\end{exa}

\begin{exa}
We compare the flip signatures of $(A, J)$, $(B, I)$ and $(B, K)$.	
Direct computation shows that the index sets for the eventual kernels of $A$ and $B$ are
$$ \mathcal{I}nd(\mathcal{K}(A)) = \{1, 6\}, \qquad \text{and} \qquad  \mathcal{I}nd(\mathcal{K}(B)) = \{1\}$$
and the flip signatures are
$$\emph{F.Sig}(A, J) = (-1, +1),$$ $$\emph{F.Sig}(B, I) = (+1)$$ 
and 
$$\emph{F.Sig}(B, K) = (-1).$$
By Theorem D, we see that
$$(\textsf{X}_A, \sigma_A, \varphi_{A, J}) \ncong (\text{X}_B, \sigma_B, \varphi_{B, I}),$$
$$(\textsf{X}_A, \sigma_A, \varphi_{A, J}) \ncong (\text{X}_B, \sigma_B, \varphi_{B, K})$$
and
$$(\textsf{X}_A, \sigma_A, \varphi_{B, I}) \ncong (\text{X}_B, \sigma_B, \varphi_{B, K}).$$
\end{exa}

In the following example, we see that the coincidence of the Lind zeta functions does not guarantee the existence of $D_{\infty}$-SE between the corresponding flip pairs.
\begin{exa}
Let
$$A=\left[\begin{array}{rrrrrrr} 1 & 1& 1 & 0 & 0 & 0 & 0 \\ 0 & 1 & 0 & 1 & 0 & 0 & 0 \\ 0 & 0 & 1 & 0 & 0 & 1 & 0 \\ 0 & 0 & 0 & 1 &0 & 0 & 1 \\ 1 & 1 & 1 & 0 & 1 & 0 & 0 \\ 1 & 1 & 1 & 0 & 0 & 1 & 0 \\ 0 & 0 & 0 & 1 & 1 & 0 & 1\end{array}\right], \quad B=\left[\begin{array}{rrrrrrr} 1 & 1& 0 & 0 & 0 & 0 & 0 \\ 0 & 1 & 0 & 1 & 1 & 1 & 0 \\ 0 & 0 & 1 & 1 & 1 & 1 & 0 \\ 0 & 0 & 0 & 1 &0 & 0 & 1 \\ 1 & 0 & 0 & 0 & 1 & 0 & 0 \\ 0 & 0 & 1 & 0 & 0 & 1 & 0 \\ 0 & 0 & 0 & 1 & 1 & 1 & 1\end{array}\right]$$
and 
$$J=\left[\begin{array}{rrrrrrr} 1 & 0 & 0 & 0 & 0 & 0 & 0 \\ 0 & 0 & 0 & 0 & 1 & 0 & 0 \\ 0 & 0 & 0 & 0 & 0 & 1 & 0 \\ 0 & 0 & 0 & 0 & 0 & 0 & 1 \\ 0 & 1 & 0 & 0 & 0 & 0 & 0 \\ 0 & 0 & 1 & 0 & 0 & 0 & 0 \\ 0 & 0 & 0 & 1 & 0 & 0 & 0 \end{array}\right].$$
The characteristic functions $\chi_A$ and $\chi_B$ of $A$ and $B$ are the same:
$$\chi_A(t) = \chi_B(t)=t(t-1)^4(t^2-3t+1).$$
We denote the zeros of $t^2 - 3t +1$ by $\lambda$ and $\mu$.
Direct computation shows that $(A, J)$ and $(B, J)$ are flip pairs and
$(\textsf{X}_A, \sigma_A, \varphi_{A, J})$ and $(\textsf{X}_B, \sigma_B, \varphi_{B, J})$
share the same numbers of fixed points. 
\begin{eqnarray*}
	p_m &=& 4 + \lambda^m + \mu^m,\\
	p_{2m-1, 0} &=& \frac{8 \lambda^{m}-3\lambda^{m-1}}{11 \lambda -4} + \frac{8 \mu^{m}-3 \mu^{m-1}}{11 \mu -4},\\
	p_{2m, 0} &=& \frac{\lambda^{m+1}}{11 \lambda -4} + \frac{\mu^{m+1}}{11 \mu -4},\\
	p_{2m, 1} &=& \frac{55\lambda^m -21\lambda^{m-1}}{11 \lambda -4} + \frac{55\mu^m -21\mu^{m-1}}{11 \mu -4} \qquad (m=1, 2, \cdots).
\end{eqnarray*} 
As a result, they share the same Lind zeta functions:
$$\sqrt{\frac{1}{t^2(1-t^2)^4(1-3t^2+t^4)}} \exp \Big(\frac{t+3t^2-t^3-2t^4}{1-3t^2+t^4}\Big).$$

If there is a $D_{\infty}$-SE $(D, E)$ from $(A, J)$ to $(B, J)$, then $(D, E)$ also becomes a SE from $A$ to $B$. 
It is well known \cite{LM} that the existence of SE from $A$ to $B$ implies that $A$ and $B$ 
have the same Jordan forms away from zero up to the order of Jordan blocks.
The Jordan canonical forms of $A$ and $B$ are given by
$$\left[\begin{array}{rrrrrrr} \lambda & & & & & & \\ & \mu & & & & & \\ & & 1 & 1 & 0 & 0 & \\ & & 0 & 1 &1 & 0 & \\  & & 0 & 0 & 1 &1 & \\ & & 0 & 0 & 0 & 1 & \\ & & & & & & 0 \end{array}\right] \quad \text{and} \quad \left[\begin{array}{rrrrrrr} \lambda & & & & & & \\ & \mu & & & & & \\ & & 1 & 1 & & & \\ & & 0 & 1 & & & \\  & & & & 1 &1 & \\ & & & & 0 & 1 & \\ & & & & & & 0  \end{array}\right],$$  
respectively.
From this, we see that $(A, J)$ cannot be $D_{\infty}$-shift equivalent to $(B, J)$.
\end{exa}

\end{document}